\documentclass[12pt,oneside]{book}
\usepackage{graphicx}
\usepackage{amsfonts}
\usepackage{amscd}
\usepackage{amssymb}
\usepackage{xy}
\usepackage{tikz}
\usepackage{lscape}
\usepackage{setspace}
\usepackage{amsthm}
\usepackage{mathtools}
\usetikzlibrary{matrix,arrows}
\newtheorem{thm}{Theorem}[section]
\newtheorem{cor}[thm]{Corollary}
\newtheorem{lem}[thm]{Lemma}
\newtheorem{prop}[thm]{Proposition}
\newtheorem{defn}[thm]{Definition}
\newtheorem{rem}[thm]{Remark}

\newtheorem{ex}[thm]{Example}
\newtheorem{conj}[thm]{Conjecture}

\newcommand{\NN}{\mathbb N}
\newcommand{\ZZ}{\mathbb Z}
\newcommand{\QQ}{\mathbb Q}

\newcommand{\CC}{\mathbb C}

\begin{document}
\frontmatter
\begin{titlepage}
    \begin{center}
        UNIVERSITY OF CALIFORNIA\\
        \vspace*{-.4cm}
        SANTA CRUZ
        
        \textbf{STRING HOMOLOGY AND LIE ALGEBRA STRUCTURES}

        A dissertation submitted in partial satisfaction
of the\\ 
\vspace*{-.4cm}
requirements for the degree of\\
DOCTOR OF PHILOSOPHY\\
in\\
MATHEMATICS\\
by\\
\textbf{Felicia Y. Tabing}

        June 2015
        \vspace*{.7cm}
    \end{center}
   \begin{flushright}
   The Dissertation of Felicia Y. Tabing\\
\vspace*{-.4cm}
is approved:\hspace*{4cm}~~~~~ \\
\vspace*{.3cm}
\noindent\rule{6.1cm}{0.4pt}\\
\vspace*{-.4cm}
 Professor Hirotaka Tamanoi~~~~~~~~\\
\vspace{.5cm}
\noindent\rule{6.1cm}{0.4pt}\\
\vspace*{-.4cm}
Professor Geoffrey Mason~~~~~~~~~~~\\
\vspace{.5cm}
\noindent\rule{6.1cm}{0.4pt}\\
\vspace*{-.4cm}
Professor Richard Montgomery~~\hspace*{1mm}
   \end{flushright}
   \begin{flushleft}
   \vspace{.8cm}
\noindent\rule{6cm}{0.4pt}\\
\vspace*{-.4cm}
Tyrus Miller\\
\vspace*{-.4cm}
Vice Provost and Dean of Graduate Studies
   \end{flushleft}
\end{titlepage}

\newpage

\thispagestyle{plain}
\pagenumbering{gobble}
\thispagestyle{empty}
\vspace*{\fill}
\begin{center}
\Large
Copyright $\textcopyright$ by \\
\vspace*{1 cm}
Felicia Y. Tabing\\
\vspace*{1 cm}
2015
\vspace*{\fill}
\end{center}

\pagenumbering{roman}
\setcounter{page}{2}
{\normalsize
\tableofcontents}
\listoffigures

\addcontentsline{toc}{chapter}{Abtract}
\newpage
\thispagestyle{plain}
\begin{center}
\Large
\textbf{Abstract}\\
\vspace{0.4cm}
    \large
    String Homology and Lie Algebra Structures
    
    \vspace{0.4cm}
    by\\
     \vspace{0.4cm}
    \textbf{Felicia Y. Tabing}

\end{center}
~
\indent Chas and Sullivan introduced string homology in \cite{CS}, which is the equivariant homology of the loop space with the $S^1$ action on loops by rotation. Craig Westerland computed the string homology for spheres with coefficients in $\ZZ /2\ZZ$ \cite{We} and in Somnath Basu's dissertation \cite{Ba}, he computes the string homology and string bracket for spheres over rational coefficients, and he finds that the bracket is trivial. In this paper, we compute string homology and the string bracket for spheres with integer coefficients, treating the odd- and even-dimensional cases separately. We use the Gysin sequence and Leray-Serre spectral sequence to aid in our computations. We find that over the integers, the string Lie algebra bracket structure is more interesting, and not always zero as in \cite{Ba}. The string bracket turns out to be non-zero on torsion coming from string homology. \\
\indent We also make some computations of the Goldman Lie algebra structure, and more generally, the string Lie algebra structure of closed, orientable surfaces.
\newpage
\thispagestyle{plain}
\addcontentsline{toc}{chapter}{Dedication}
\begin{center}
\vspace*{\fill}
To Michael Kusuda.
\vspace*{\fill}
\end{center}

\newpage
\thispagestyle{plain}
\addcontentsline{toc}{chapter}{Acknowledgements
}
\Large
\begin{center}
\textbf{Acknowledgements}
\end{center}
\normalsize
\vspace{0.4cm}
~\\
\indent I would like to express my deepest gratitude to my advisor Hirotaka Tamanoi for his support of me over the past few years. He has been incredibly patient with me, and I am very grateful for him pushing me to learn how to work on my own. I enjoyed the time spent in his office, learning about algebraic topology. I also very much valued hearing his views on life, which I recall when I am having a hard time. \\
\indent I wish to thank the rest of my thesis committee, Geoff Mason and Richard Montgomery. Richard, for being available to chat about mathematics, and introducing me to Bill Goldman.\\
\indent I am most grateful to Debra Lewis for the help and guidance she has given me. Her support of me through all my years as a graduate student, and even as an undergraduate, has been invaluable. She was always available to talk when I was feeling anxious and needed moral support, and provided me with encouragement.\\
\indent Teaching has been one the most enjoyable part of graduate school. I would like to acknowledge Frank Ba{\"u}erle, who was a  teaching mentor to me. I hope one day that I can be as much of a compassionate and patient teacher as he is.\\
\indent  I would like to thank the following fellow graduate students, current and former, for their friendship, support, mathematical and non-mathematical conversations: Alex Beloi, Victor Bermudez, Jonathan Chi, Michael Campbell, Jamison Barsotti, Sean Gasiorek, Rob Carman, Gabriel Martins, Mitchell Owen, Vinod Sastry, Shawn Tsosie, and Wei Yuan. Liz Pannell deserves a special mention for her friendship and starting the Noetherian Ring with me. I am grateful to Danquynh Nguyen, for being so generous and sharing her delicious food with me, as I would often come to campus without food, and very hungry. Jean Verrette, whom I am grateful to have been paired with as roommates at the Algebraic Topology Summer School at MSRI.\\
\indent I am grateful for the support of my parents, Sylvia and German Tabing, for not discouraging me from mathematics, even though when I was a kid I said I wanted to be surgeon. I want to acknowledge Linda and Harry Kusuda, for their support and believing in me. I also want to thank Guy Gov and Annie Nguyen, with whom I can forget about my academic worries and have fun.\\
\indent Lastly, I am greatly indebted to Michael Kusuda for his never-ending love an support. He was my greatest supporter, made sure I was well fed, and took care of my every need. With his support, I was allowed to concentrate on learning mathematics, and I am extremely grateful to him.
\mainmatter
\setcounter{chapter}{-1}
\chapter{Introduction}
~
\indent The term \emph{String Topology} came from the paper of the same name by Moira Chas and Dennis Sullivan in 1999. This paper discussed the various algebraic structures that arose from the homology of the free loop space. This paper came out of trying to generalize the Lie algebra structure that William M. Goldman described by the intersection and concatenation of loops on surfaces \cite{Go}.\\
\indent Chapter 1 is an introduction to the Goldman Lie algebra, and we explore its structure. In particular we consider the structure of the Lie algebra for the closed torus, including computations showing it is finitely generated. \\
\indent Chapter 2 introduces string topology background needed for the rest of this paper, and the various algebra structures of loop homology and string homology.\\
\indent Chapter 3 contains the computations of the integral string homology and string bracket structure for spheres, where some torsion phenomena appear. In our computations, we use the Leray-Serre spectral sequence, and the Gysin exact sequence.\\
\indent Chapter 4 explores the string homology and bracket structure of surfaces. 
\chapter{The Goldman Lie Algebra}
~
\indent The Goldman Lie algebra was introduced by William M. Goldman in 1986 \cite{Go}. \\
\indent Throughout, let $\Sigma_{g,n}$ denote an oriented, genus $g$ surface with $n \geq 0$ boundary components. Denote $\hat{\pi}(\Sigma_{g,n})$ to be the set of free homotopy classes of loops on $\Sigma_{g,n}$, where the surface is not mentioned in the notation of $\hat{\pi}$ when it is clear from the context that we are talking about some fixed surface. Recall the following.
\lem The set of free homotopy classes of loops on a surface $\Sigma_{g,n}$ is in one-to-one correspondence with conjugacy classes of $\pi_1(\Sigma_{g,n})$.
\rem \normalfont We can represent homotopy classes of loops by cyclically reduced words with letters the generators of the fundamental group.
\defn \normalfont Fix a surface $\Sigma_{g,n}$ and an orientation of $\Sigma_{g,n}$. Let $\alpha, \beta \in \hat{\pi}$. The \textbf{Goldman bracket} of $\alpha$ and $\beta$ is defined to be
\begin{align}
[\alpha, \beta]=\sum_{p \in \alpha \cap \beta} \epsilon (p) \alpha *_p \beta
\end{align}
where $\alpha$ and $\beta$ intersect in transverse double points $p$, and $\epsilon (p)$ is the sign of the intersection, or $\epsilon (p)=1$ if the ordered vectors in the tangent space to $\Sigma_{g,n}$ tangent to loop $\alpha$ and $\beta$ match the orientation of the surface, and $\epsilon (p)=-1$ otherwise. 
\begin{figure}
\includegraphics[scale=.4]{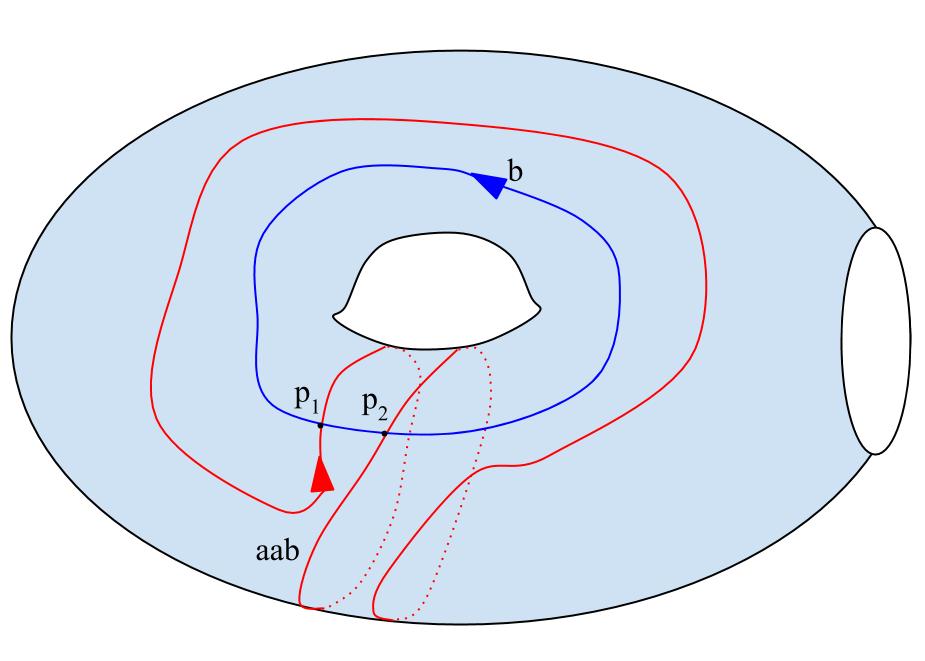}
\caption{Loops $aab$ and $b$ on the torus with one boundary component.\label{fig:sigma11}}
\end{figure}
\ex \normalfont We compute $[aab,b]$ on the surface $\Sigma_{1,1}$. The loops represented by words $aab$ and $b$ are shown in Figure~\ref{fig:sigma11}, with intersection points $p_1$ and $p_2$. At the intersection point $p_1$, we smooth the intersection by creating a new loop, $aabb$, by following the red loop $aab$ in the direction of its orientation at $p_1$, and when returning to $p_1$, we now follow the blue loop $b$ in the direction of its orientation. When we return back to $p_1$, we close the loop. At the intersection $p_2$, we do the same, and create the loop $abab$.  We get that $[aab,b]=\pm(aabb+abab)$ where the sign depends on the chosen orientation of $\Sigma_{1,1}$.
\thm (Goldman) The Goldman bracket is well defined, skew-symmetric, and satisfies the Jacobi identity
\proof \cite{Go} \qed\\
\indent We can extend the bracket linearly to $\ZZ [\hat{\pi}]$ (or $\QQ [\hat{\pi}]$), the free module over $\ZZ$ (or $\QQ$) with basis $\hat{\pi}$, to get a bilinear map
\begin{align*}
[ -,-]: \ZZ \hat{\pi}\times \ZZ \hat{\pi} \rightarrow \ZZ \hat{\pi}
\end{align*}.
Thus, $\ZZ \hat{\pi}$ is a Lie algebra with bracket $[-,-]$, which we call the \emph{Goldman Lie Algebra}, denoted by $\mathfrak{G}$ throughout the rest of this chapter. When it is unclear what the surface we are referring to, we use $\mathfrak{G}_{\Sigma_{g,n}}$
\section{Goldman Lie Algebra Structure}
~
\indent What is the Lie Algebra structure of the Goldman Lie algebra? So far, the center of the Goldman Lie algebra is known, but much of the structure is still a mystery.
\thm (Etingof) The center of $\mathfrak{G}_{\Sigma_{g,0}}$ is spanned by the contractible loop \cite{Et}.
\thm (Kabiraj) The center of $\mathfrak{G}_{\Sigma_{g,n}}$ is generated by peripheral loops \cite{Ka}\\
\normalfont
\indent A question posed by Chas \cite{Ch} is whether or not $\mathfrak{G}$ is finitely generated. In Goldman's paper \cite{Go}, he also introduces what is called the homological Goldman Lie algebra. This Lie algebra is defined on intersection form on the first homology group of a surface. It is known that this Lie algebra is indeed finitely generated \cite{KKT}, but of course, the homological Goldman Lie algebra is simpler. \\
\indent The closed torus is a special case. $\mathfrak{G}_{\Sigma_{1,0}}$ is finitely generated. Recall that we can represent free homotopy classes of loops on $\Sigma_{1,0}$ by cyclically reduced words in two letters, $a$ and $b$, and we can represent all homotopy classes of loops on the torus by the word $a^lb^k$ for $k,l\in \ZZ$.
\prop \label{prop:torusstructure} The Goldman bracket structure of $\mathfrak{G}_{\Sigma_{0,1}}$ is given by
\begin{align*}
[a^ib^j,a^kb^l]=(il-jk)a^{i+k}b^{j+l}
\end{align*}
\thm \label{thm:torus}$\mathfrak{G}_{\Sigma_{1,0}}$ is finitely generated when considered as a Lie algebra over $\QQ$..
\proof  We denote a contractible loop by $1$. We claim that $\mathfrak{G}_{\Sigma_{1,0}}$ is generated by $\{a,b,a^{-1},b^{-1}\}$. This will take many steps. We will first show that we can generate certain homotopy classes of loops. Below, we assume $n\neq 0$.
\begin{enumerate}
\item $a^nb^1=[a,a^{n-1}b]$, which we get inductively,
\item $a^n=[b^{-1},-\frac{1}{n}a^nb]$ 
\item $ab^n=[b,-ab^{n-1}]$
\item $b^n=[a^{-1},-\frac{1}{n}ab^n]$
\item $a^nb^n=[a^n,\frac{1}{n^2}b^n]$
\item $a^{-n}b=[a^{-1},-a^{-n+1}b]$
\item $a^{-n}=[b^{-1},\frac{1}{n}a^{-n}b]$
\item $a^{-1}b^n-[a^{-1},-\frac{1}{n}b^n]$
\item $a^{-n}b^n=[a^{-n},-\frac{1}{n^2}b^n]$
\item $ab^{-n}=[b^{-1},ab^{n+1}]$ which we get inductively,
\item $b^{-1}=[a^{-1},\frac{1}{n}ab^{-n}]$
\item $a^{-n}b^{-n}=[a^{-n},\frac{1}{n^2}b^{-n}]$
\item $a^nb^{-n}=[a^n,-\frac{1}{n^2}b^{-n}]$.
\item From 13. and 9. for $n=1$, we get $a^0b^0=1=[ab^{-1},\frac{1}{2}a^{-1}b]$.
\end{enumerate}
We still have a few more cases to show, namely how to generate the homotopy class of the loop $a^ib^j$ in the following cases.
\begin{enumerate}
\item[Case 1:] Suppose $i,j>0$.
\begin{enumerate}
\item Suppose $i<j$, then $j=i+r$ for some $r\in \ZZ-\{0\}$.\\ Then $\frac{1}{ar}[a^ib^i,b^j]=a^ib^j$.
\item Suppose $i>j$, then $i=j+r$ for $r\in \ZZ-\{0\}$.\\ Then $-\frac{1}{br}[a^jb^j,a^r]=a^ib^j$.
\end{enumerate}
\item[Case 2:] Suppose $i<0<j$.
\begin{enumerate}
\item Suppose $|i|<|j|$, then $j=-i+r$ for $r\in \ZZ-\{0\}$.\\ Then $\frac{1}{ar}[a^ib^{-i},b^r]=a^ib^j$.
\item Suppose $|a|>|b|$, then $i=-j+r$ for $r\in \ZZ-\{0\}$.\\ Then $-\frac{1}{br}[a^{-j}b^{j},a^r]=a^ib^j$.
\end{enumerate}
\item[Case 3:] The case $i,j<0$, and $i\neq j$ is similar to Case 1.
\item[Case 4:] The case $b<0<a$ is similar to Case 2.
\end{enumerate}
Thus, everything in $\mathfrak{G}_{\Sigma_{1,0}}$ can be generated as a Lie algebra over $\QQ$.
\qed
\cor We can refine the generators of $\mathfrak{G}_{\Sigma_{0,1}}$ to a smaller basis, namely $\{a, a^{-1}b^{-1}+b+1,b\}$
\proof We show that we generate the basis elements mentioned in the proof of Theorem~\ref{thm:torus}.
\begin{align*}
[a^{-1}b^{-1}+b+1,b]&=-a^{-1},\\
[a^{-1}b^{-1}+b+1,a]&=b^{-1}-ab,\\
[a,b]&=ab.
\end{align*}
\qed
\cor $\mathfrak{G}_{\Sigma_{1,0}}$ as a Lie algebra over $\QQ$ is not nilpotent, nor solvable, since $[\mathfrak{G}_{\Sigma_{1,0}} , \mathfrak{G}_{\Sigma_{1,0}}] =\mathfrak{G}_{\Sigma_{1,0}}$.
\rem $\mathfrak{G}_{\Sigma_{1,0}}$ is not finitely generated as a Lie algebra over $\mathbb{Z}$.
\proof We will show that the set $\{(n-1)a^n\}_{n>2, n\in \ZZ}$ cannot be generated. Suppose to the contrary that we can generate $(n-1)a^n$, so there exists $(i_s,j_s),(k_s,l_s) \in \ZZ ^2$ such that
\begin{align*}
\sum_{s=1}^t\pm [a^{i_s}b^{j_s},a^{k_s},b^{l_s}]=(n-1)a^n.
\end{align*}
Then, as in Proposition~\ref{prop:torusstructure},
\begin{align*}
\sum_{s=1}^t\pm [a^{i_s}b^{j_s},a^{k_s},b^{l_s}]=\sum_{s=1}^t\pm (i_sl_s-j_sk_s)a^{i_s+k_s}b^{j_s+l_s}.
\end{align*}
We need that $i_s+k_s=n$ and $j_s+l_s=0$, so $i_sl_s-j_sk_s=-i_sj_s-j_sn+j_si_s=-j_sn$. So $n\mid i_sl_s-j_sk_s$, and $n\mid \sum_{s=1}^t\pm (i_sl_s-j_sk_s)$, so $n\mid (n-1)$, which is a contradiction. \qed
\conj We conjecture that $\mathfrak{G}_{\Sigma_{g,n}}$ for $g\geq 1$ and $n>1$ is not finitely generated. For the particular case for a punctured torus, the peripheral loop is given by a commutator word. We noticed in using Chas' program for computing the bracket seems to not generate a commutator word, nor products of commutators. This needs more work, but this would mean we have a set $\{(aba^{-1}b^{-1})^n\}_{n \in \ZZ}$ of infinitely many homotopy classes of loops that each cannot be generated by any other homotopy classes of loops. 
\normalfont
\prop \label{prop:derived} The derived Lie algebra for $\mathfrak{G}_{\Sigma_{1,0}}$ is given by
\begin{align*}
[\mathfrak{G}_{\Sigma_{1,0}},\mathfrak{G}_{\Sigma_{1,0}}]=\langle d(a^ib^j),na^n,nb^n \rangle
\end{align*}
for $d=gcd(i,j)$ and $n\in \ZZ-\{0\}$.
\proof We first show $[\mathfrak{G}_{\Sigma_{1,0}},\mathfrak{G}_{\Sigma_{1,0}}]\subset \langle d(a^ib^j),na^n,nb^n \rangle$
\begin{enumerate}
\item[Case 1:] Suppose $d=gcd(i,j)$, $i,j \neq 0$ and $ma^ib^j \in [\mathfrak{G}_{\Sigma_{1,0}},\mathfrak{G}_{\Sigma_{1,0}}]$ for some $m\in \ZZ$. Write $xi+yj=d$ for some $x,y\in \ZZ$ and $ma^ib^j=[a^kb^l,a^p,b^q]$ for $k,l,p,q \in \ZZ$. But 
\begin{align} \label{align:gcdd}
[a^kb^l,a^p,b^q]=(kq-lp)a^{k+p}b^{l+q}
\end{align}
so we have that $k+p=i$, $l+q=j$, and $kq-lp=kj-li=d(k(\frac{j}{d})-l(\frac{i}{d}))=m$. Thus $d\mid m$.
\item[Case 2:] Suppose that $n \neq 0$ and that $[a^ib^j,a^kb^l]=ma^n$ for $i,j,k,l,m\in \ZZ$. Then $i+k=n$, $j+l=0$, so
\begin{align} \label{align:n0}
[a^ib^j,a^kb^l]=-jna^n
\end{align} 
so $n\mid m$. 
\item[Case 3:] Showing that for $n \neq 0$ and $n\mid m$ for $mb^n \in [\mathfrak{G}_{\Sigma_{1,0}},\mathfrak{G}_{\Sigma_{1,0}}]$ is similar to Case 2.\\
\indent To show the other containment, we can consider the equality~\ref{align:gcdd} with $k=y$ and $l=-x$ for Case 1, we can consider the equality~\ref{align:n0} with $j=-1$, and we can do something similar for Case 3.
\end{enumerate}
\qed
\prop The lower central series for $\mathfrak{G}_{\Sigma_{1,0}}$ stabilizes, i.e.
\begin{align*}
[\mathfrak{G}_{\Sigma_{1,0}},G_i]=\langle d(a^ib^j),na^n,nb^n \rangle
\end{align*}
where $d=gcd(i,j)$, $n\in \ZZ-\{0\}$, for all $i\geq 0$, and $G_i=[\mathfrak{G}_{\Sigma_{1,0}},G_{i-1}]$ defined inductively, where $G_0=\mathfrak{G}_{\Sigma_{1,0}}$.
\proof For $i=1$, this is just Proposition~\ref{prop:derived}. For $i=2$, we need to show that 
\begin{align*}
[ \mathfrak{G}_{\Sigma_{1,0}}, [\mathfrak{G}_{\Sigma_{1,0}}, \mathfrak{G}_{\Sigma_{1,0}} ]]=\langle d(a^ib^j),na^n,nb^n \rangle .
\end{align*}
The "$\subset $" containment is clear. First, consider $a^ib^{-1}\in \mathfrak{G}_{\Sigma_{1,0}}$ and $a^{n-i}b \in \langle d(a^ib^j),na^n,nb^n\rangle$ (since $gcd(n-i,1)=1$). We have that 
\begin{align*}
[a^ib^{-1},a^{n-i}b]=na^n.
\end{align*}
In a similar way, we can show that $nb^n \in [ \mathfrak{G}_{\Sigma_{1,0}}, [\mathfrak{G}_{\Sigma_{1,0}}, \mathfrak{G}_{\Sigma_{1,0}} ]]$. \\
Now consider $d=gcd(i.j)$, so we can write $d=xi+yj$. Consider $a^{i+y}b^{j-x} \in \mathfrak{G}_{\Sigma_{1,0}}$ and $ a^{-y}b^x \in [\mathfrak{G}_{\Sigma_{1,0}} , \mathfrak{G}_{\Sigma_{1,0}} ]$ (since $1=\frac{i}{d}x+\frac{j}{d}y$ implies $gcd(x,y)=1$. We have
\begin{align*}
[a^{i+y}b^{j-x} , a^{-y}b^x]=da^ib^j.
\end{align*}
Thus, it follows that the lower central series stabilizes.
\qed
\cor $\mathfrak{G}_{\Sigma_{1,0}}$ as a Lie algebra over $\ZZ$ is not nilpotent.
\normalfont
\chapter{String Topology Preliminaries}
~
\indent Here we describe the basic algebraic structures appearing in the homology and equivariant homology of the free loop space, as described by Chas and Sullivan in \emph{String Topology}. Throughout this paper, let $M$ be a manifold of dimension $d$, $\Omega M$ the based loop space of $M$, and the space of all continuous, piecewise smooth loops on $M$,  $LM=Map(S^1, M)$, the free loop space of $M$. Note that $LM$ can be considered to be an infinite-dimensional manifold, and it is topologised with the compact-open topology. We will consider homology and cohomology with integer coefficients, unless otherwise stated. We denote the usual homology of the free loop space of $M$ as $H_*(LM)$ and equivariant homology will be denoted by $H_*^{S^1}(LM)$. 
\section{Loop Homology Algebra Preliminaries}
~
\indent We first describe the Chas-Sullivan loop product, which Chas and Sullivan defined on the chain level of $LM$, the space of all continuous, piecewise smooth loops on $M$. \\
\indent The loop product is a combination of the intersection product and the product given by the concatenation of loops. It is defined transversally at the chain level.\\
\indent Consider an $i$-chain of loops in $LM$. We can think of a simplex in this chain as a map $\sigma : \Delta_i \rightarrow LM$ or as a map $\sigma : \Delta_i \otimes S^1 \rightarrow M$. So we can think of an $i$ chain of loops as a map from a simplex with loops above it into $M$. Intuitively, if we have an $i$-chain and a $j$-chain of loops where the marked points intersect transversally, then we get a new $i+j-d$-chain of loops consisting of the intersecting marked points, and at each marked point, the new loop is formed by going around the $i$-chain loops then around the $j$-chain of loops. This description at the chain level can pass to homology to form the \emph{Chas-Sullivan loop product}.\\
\indent Here we give a more precise description of the product given in \emph{String Topology and Cyclic Homology} \cite{CHV}. Let $Map(8,M)=Map(S^1\vee S^1)$ be the space of continuous, piecewise smooth maps from the figure eight, or the wedge sum of two circles to $M$. This is topologised with the compact-open topology and can be considered as an infinite-dimensional manifold, but we need piecewise smooth in order for it to be some sort of manifold. It can also be viewed as a subspace of $LM \times LM$ where the loops agree at $0$.\\
\indent Consider the following diagram. The left square is a pullback diagram.:

\begin{equation*}\begin{CD}
LM \times LM @< \rho_{\text{in}}<< Map(8,M) @> \rho_{\text{out}} >> LM\\
@VV{ev \times ev}V @VV{ev}V\\
M \times M @<\Delta << M
\end{CD}\end{equation*}
\\where $\rho_{\text{in}}$ is the restriction of the figure eight to the two different loops, and $\rho_{\text{out}}$ is where the figure eight loop is considered as one loop. The maps $ev$ are the evaluation of the loops at $0$, and $\Delta$ is the diagonal map. $ev$ is a locally trivial fibration so if $\eta_{\Delta}$ is a tubular neighborhood of the diagonal embedding, then a tubular neighborhood of the $\eta_{\rho_{in}}=(ev\times ev)^{-1}(\eta_{\Delta})$ is homeomorphic to $ev^*(TM)=ev^*(\Delta(M))$. We can have a tubular neighborhood since $\rho_{in}$ is a codimension $d$ embedding. Since $ev^*(\eta_{\Delta})$ is the pullback of $\eta_{\Delta}$, which has fiber dimension $d$ since it is the normal bundle, the pullback $ev^*(\eta_{\Delta})$ has fibers isomorphic to fibers of $\eta_{\Delta}$, so $ev^*(\eta_{\Delta})$ also has fiber dimension $d$. This means that the normal bundle of $map(8,M)$ in $LM\times LM$ has codimension $d$. Since in the above diagram, the left square is a pullback diagram of fiber bundles, we have that $\rho_{in}$ is a codimension $d$ embedding.\\
\indent The induced maps on homology go in the wrong direction, in order to remedy this, we need to turn the map $\rho_{\text{in}}$ around somehow. We do this by using the Pontrjagin-Thom collapse map:\\
\begin{equation*}
LM \times LM \rightarrow LM \times LM / LM \times LM -ev^*(TM) \cong Thom(Map(8,M))
\end{equation*}\\
Define the umkehr map $(\rho_{in})_!$ containing the induced map on homology above as follows:\\
\begin{align*}
(\rho_{in})_!: H_*(LM)\otimes H_*(LM)\cong H_*(LM \times LM)\rightarrow & H_*(Thom(Map(8,M)))\\ &\cong H_{*-d}(Map(8,M))
\end{align*}
\indent The last isomorphism is given by the Thom Isomorphism by taking the cap product with $u\in H^d(Thom(Map(8,M)))$, the Thom class given by the orientation.
\defn \normalfont \label{def:loopproduct} The following composition gives the \emph{Chas-Sullivan loop product} (or just \emph{loop product}):
\begin{align}
-\bullet - =(\rho_{out})_* \circ (\rho_{in})_!: H_*(LM \times LM) \rightarrow H_{*-d}(Map(8,M)) \rightarrow H_{*-d}(LM)
\end{align}
This product can be extended to homology. It is convenient to regrade the loop homology as follows:
\begin{align*}
\mathbb{H}_*(LM):=H_{*+d}(LM)
\end{align*}
\\we can rewrite the product:
\begin{equation*}
-\bullet -= \mathbb{H}_*(LM) \otimes\mathbb{H}_*(LM) \rightarrow \mathbb{H}_{*}(LM)
\end{equation*}
We may drop the $LM$ from the notation and denote loop homology by $\mathbb{H}_*$ when it is clear which manifold we are referring to.
\thm(Chas-Sullivan) $(\mathbb{H}_*(LM),\bullet )$ is an associative, graded, commutative algebra.
\normalfont
\defn \normalfont There is a \emph{Batalin-Vilkovisky operator} denoted by $\Delta$, which comes from the natural action given by rotation of loops,
\begin{align*}
\rho:S^1\times LM \rightarrow LM
\end{align*} 
given by $\rho(t,\gamma)(s)=\gamma(s+t)$. This action defines a degree one operator on loop homology: $\Delta:\mathbb{H}_*(LM)\rightarrow \mathbb{H}_{*+1}(LM)$ given by $\delta(\alpha)=\rho_*([S^1]\otimes \alpha)$ for $\alpha \in H_k(LM)$. 
\thm (Chas-Sullivan) \label{thm:loopbracket} $(\mathbb{H}_*(LM), \Delta )$ is a Batalin-Vilkovisky algebra, 
\begin{enumerate}
\item $(\mathbb{H}_*(LM), \bullet)$ is a graded, commutative, associative algebra
\item $\Delta \circ \Delta =0$
\item $(-1)^{|\alpha |} \Delta (\alpha \bullet \beta ) - (-1)^{ |\alpha |} \Delta (\alpha) \bullet \beta - \alpha \bullet \Delta (\beta)$ is a derivation in each variable.
\end{enumerate}
\normalfont
We can also define a Lie bracket with $\bullet$ and $\Delta$ as in part 3 of Theorem~\ref{thm:loopbracket}.
\defn \normalfont The \emph{loop bracket} is defined as
\begin{align*}
\{\alpha , \beta \}:= (-1)^{|\alpha |} \Delta (\alpha \bullet \beta ) - (-1)^{ |\alpha |} \Delta (\alpha) \bullet \beta - \alpha \bullet \Delta (\beta)
\end{align*}
which is the deviation of $\Delta$ from being a derivation of $\bullet$.
\thm (Chas-Sullivan) $(\mathbb{H}_*(LM),\bullet , \{-,-\})$ has the structure of a Gerstenhaber algebra,
\begin{enumerate}
\item $(\mathbb{H}_*(LM), \bullet)$ is a graded, commutative, associative algebra
\item $\{-,-\}$ is a degree $1$ Lie bracket,
\begin{enumerate}
\item ${\alpha , \beta }=(-1)^{ (|\alpha | +1)(|\beta |+1)+1} \{\beta ,\alpha \}$,
\item $\{ \alpha , \{ \beta , \gamma \} \} =\{ \{ \alpha, \beta \} , \gamma \} + (-1)^{ (|\alpha | +1)(|\beta |+1)} \{ \beta , \{ \alpha , \gamma \} \}$,
\end{enumerate}
\item $\{ \alpha , \beta \bullet \gamma \} = \{\alpha , \beta \} \bullet \gamma + (-1)^{ (|\alpha | -1)|\beta |} \beta \bullet \{ \alpha , \gamma \}$.
\end{enumerate}
\normalfont
\section{String Homology Algebra Preliminaries}
~
\indent Now we consider algebraic structures on the equivariant homology of the free loop space with respect to the action of rotation of loops, $H_*^{S^1}(LM)$. Consider the fibration
\begin{align*}
S^2 \rightarrow LM \times ES^1 \rightarrow LM \times_{S^1} ES^1.
\end{align*}
This induces a long exact sequence on homology, the Gysin sequence from which we will use to describe a Lie bracket on $H_*^{S^1}(LM)$.
\begin{align*}
\cdots \rightarrow \mathbb{H}_{*-d}(LM) \xrightarrow{e} H_*^{S^1}(LM) \xrightarrow{\cap} H_{*-2}^{S^1}(LM) \xrightarrow{M} \mathbb{H}_{*-d-1}(LM) \rightarrow \cdots
\end{align*}
where $e$ and $M$ are informally called the "erasing map" and "marking map," respectively. The map $e$ forgets the marked points on the loops, and the map $M$ puts markings back on the loops in all possible places.  We have that $M$ is a homomorphism of graded Lie algebras, it preserves the brackets, going from the string bracket to the loop bracket. The map $e$ is the induced fibration map. For the rest of this paper, it will be clear from the context which space we are referring to, so we often drop the $LM$ from the homology notation.
\rem Note that $e\circ M=0$ by exactness, and $\Delta = M \circ e$.
\defn \normalfont \label{defn:stringbracket} For two classes $\alpha , \beta \in H_*^{S^1}(LM)$, we can define the \emph{string bracket} by
\begin{align*}
[\alpha , \beta ] = (-1)^{| \alpha | -d }e(M(\alpha) \bullet M(\beta) ) 
\end{align*}
where $\bullet$ was the loop product mentioned in Definition~\ref{def:loopproduct}.
\thm (Chas-Sullivan) $(H_*^{S^1}(LM),[-,-])$ is a graded Lie algebra, with Lie bracket of degree $2-d$.\\
\normalfont
More precisely, our bracket is a map:
\begin{align*}
[-,-]: H_i^{S^1}(LM)\times H_j^{S^1}(LM) \rightarrow H_{i+j+2-d}^{S^1}(LM).
\end{align*}
\indent In the following chapters, we compute the $H_*^{S^1}(LS^n)$ for all $n\in \NN$ and we compute the structure of the string bracket. 
\chapter{Computations of String Homology and the String Bracket}
~
\indent In this chapter, we compute explicitly the integral string homology and the string bracket for spheres. Somnath Basu made some computations of rational string homology for spheres in his Ph.D. thesis \cite{Ba} using rational homotopy theory and minimal models. Craig Westerland also made computations of string homology over $\ZZ_2$ for spheres in \emph{String Homology of Spheres and Projective Spaces} \cite{We} using a spectral sequence. We separate the computations for the even- and odd-dimensional spheres. First, we compute particular examples, $S^1$, $S^3$, and $S^2$, to get a better hold on the computation, then generalize to the higher-dimensional spheres. We use primarily the Gysin exact sequence, and the Leray-Serre spectral sequence to aid in computations of string homology. We find that there is a lot of interesting torsion in integral string homology, and the bracket structure is not always zero.
\section{String Homology and String Bracket of $S^1$}
~
\indent We compute the string homology of $S^1$ using the Gysin sequence for the circle bundle
\begin{eqnarray}
 \nonumber S^1 \rightarrow LS^1 \times ES^1 \rightarrow LS^1 \times_{S^1} ES^1
\end{eqnarray}\\
Basu computed this in his thesis, but here we use elementary techniques. \\
\indent Recall that the non-equivariant homology of $LS^1$ is given as follows \cite{CJY}, \cite{He}, 
\begin{equation}
\nonumber \mathbb{H}_*(LS^1)=\Lambda_{\ZZ}[a]\otimes \ZZ[x,x^{-1}], \; |a |=-1, |x|=0. 
\end{equation}
where $\mathbb{H}_*(LS^1)=H_{*+1}(LS^1)$ and $a$ corresponds to the dual of  $[S^1]$ under the geometric grading \cite{Se}, \cite{CJY}.\\
\indent The BV-operator ($\Delta=M\circ e$) acts on generators of $\mathbb{H}_*(LS^1)$ as follows, \cite{Me}: 
\begin{eqnarray}
&\Delta (a \otimes x^i)& =i(1\otimes x^{i})\nonumber \\
&\Delta (1\otimes x^i) & =0.\nonumber 
\end{eqnarray}
\indent Consider the Gysin sequence for the above circle bundle:\\
\begin{tikzpicture}[descr/.style={fill=white,inner sep=1.5pt}]
        \matrix (m) [
            matrix of math nodes,
            row sep=2em,
            column sep=2.5em,
            text height=1.5ex, text depth=0.25ex
        ]
        { &  & &0 \\
          & H^{S^1}_2(LS^1) & H^{S^1}_0(LS^1) & \mathbb{H}_{0}(LS^1)\cong \bigoplus_{n\in \ZZ} \ZZ (1\otimes x^n) \\
            & H^{S^1}_1(LS^1) & H^{S^1}_{-1}(LS^1)\cong 0 & \mathbb{H}_{-1}(LS^1)\cong \bigoplus_{n\in \ZZ} \ZZ (a\otimes x^n) \\
            & H^{S^1}_0(LS^1) & 0 & \\
           \\
        };
        \path[overlay,->, font=\scriptsize,>=latex]
    
        (m-1-4) edge[out=355,in=175,cyan] node[descr,yshift=0.3ex] {$e$} (m-2-2)
        (m-2-2) edge  node[above] {$c$} (m-2-3)
        (m-2-3) edge[orange] node[above] {$M$}(m-2-4)
        (m-2-4) edge[out=355,in=175,cyan] node[descr,yshift=0.3ex] {$e$} (m-3-2)
        (m-3-2) edge (m-3-3)
        (m-3-3) edge[orange] node[above] {$M$} (m-3-4)
        (m-3-4) edge[out=355,in=175,cyan] node[descr,yshift=0.3ex] {$e$} (m-4-2)
        (m-4-2) edge (m-4-3)
         ;
\end{tikzpicture}\\
\indent The end of the Gysin sequence gives us that ${H^{S^1}_0(LS^1)\cong \bigoplus\limits_{n\in \ZZ} \ZZ (e(a\otimes x^n))}$. Using the information from the BV-operator, ${M\circ e(a \otimes x^n)=\Delta (a \otimes x^n)}=n(1 \otimes x^n)$. Since $e$ is surjective and ${ker(e)=im(M)= \bigoplus\limits_{n\in \ZZ} n\ZZ(1 \otimes x^n)}$, we have that ${H^{S^1}_1(LS^1) \cong  \mathbb{H}_{0}(LS^1) / ker(e) \cong \bigoplus\limits_{n\in \ZZ} \ZZ /n\ZZ (1 \otimes x^n) \oplus \ZZ (1\otimes 1)} $. From the beginning of the Gysin sequence, we have ${im(c)=ker(M)=\ZZ (a\otimes 1)}$, and since $c$ is injective, ${H^{S^1}_2(LS^1) \cong \ZZ (a \otimes 1)}$. Summarizing, we get the following remark.
\rem 
\begin{eqnarray}
& H^{S^1}_0(LS^1) &\cong \bigoplus_{n\in \ZZ} \ZZ (e(a\otimes x^n))\nonumber \\
& H^{S^1}_1(LS^1) &\cong  H^{S^1}_{2i+1}(LS^1) \cong \bigoplus_{n\in \ZZ -\{0\}} \ZZ /n\ZZ (1 \otimes x^n)\oplus \ZZ (1 \otimes 1), \hspace*{.5cm} i\geq 0\nonumber \\
& H^{S^1}_2(LS^1) &\cong H^{S^1}_{2i}(LS^1) \cong \ZZ (a \otimes 1), \hspace*{.5cm} i\geq 1\nonumber 
\end{eqnarray}\\
\normalfont
\indent The string bracket, ${[-,-]: H_i^{S^1}(LS^1) \otimes H_j^{S^1}(LS^1) \rightarrow H_{i+j+1}^{S^1}(LS^1)}$ is a degree $+1$ map, and it is only nontrivial on generators of degree zero since the marking map is trivial for generators of degree greater than zero. For ${a \otimes x^n}$, ${a \otimes x^m}$ in ${H^{S^1}_0(LS^1)}$,
\begin{align*}
[e(a \otimes x^n),e( a \otimes x^m)] & =(-1)^{-1}e(M(e(a \otimes x^n))\bullet M(e(a \otimes x^m))) \\
&=-e(n(1 \otimes x^n) \bullet m(1 \otimes x^m))\\
& =-nm(e(1 \otimes x^{n+m}))\\
&=-nm(1 \otimes x^{n+m})
\end{align*}
So ${[a \otimes x^n, a \otimes x^m]=0}$ if $n+m \neq 0$ and $n+m$ divides $nm$. If $n+m=0$ then  ${[a \otimes x^n, a \otimes x^m]=nm(1\otimes 1)}$.
We can conclude that the bracket is only nontrivial for the torsion elements. 


\section{String Homology and String Bracket of $S^3$}
~
\indent We compute the equivariant homology of $S^3$ using the Gysin sequence for the circle bundle 
\begin{align}\label{eqn:threecirclebundle}
 S^1 \rightarrow LS^3 \times ES^1 \rightarrow LS^3 \times_{S^1}  ES^1
\end{align} and the Serre homology spectral sequence for 
\begin{equation}
LS^3 \rightarrow LS^3 \times_{S^1} ES^1 \rightarrow \CC P^\infty\nonumber 
\end{equation}
\indent First we compute the equivariant cohomology of $LS^3$ and then translate it to equivariant homology. We also compute the erasing ($e$) and marking ($M$) maps, as in Chas and Sullivan's paper, to compute the String Bracket.\\
\indent 
\subsection{First Few Equivariant Homology Groups of $LS^3$}
~
\indent By equivariant homology, we mean the homology of the Borel construction from the natural action of $S^1$ on $LS^3$ by rotation, denoted by $H_*^{S^1}(LS^3)=H_*(LS^3 \times_{S^1} ES^1)$. We calculate the first few equivariant homology groups of $LS^3$ to aid in our computation of the equivariant cohomology of $LS^3$.\\
\indent 
Recall that the non-equivariant homology of $LS^3$ is given as follows:
\begin{equation}
\mathbb{H}_*(LS^3)=\Lambda_{\ZZ}[\alpha]\otimes \ZZ[y], \; |\alpha |=-3, |y|=2 \nonumber 
\end{equation}
where $\mathbb{H}_*(LS^3)=H_{*+3}(LS^3)$ and $\alpha$ corresponds to the dual of  $[S^3]$ under the usual grading \cite{Se}, \cite{CJY}.\\
\indent To compute the equivariant homology of $LS^3$, we consider the Gysin sequence for the following fibration:
\begin{equation}
S^1\rightarrow LS^3 \times ES^1 \rightarrow LS^3 \times_{S^1} ES^1.\nonumber 
\end{equation}
and the BV-operator ($\Delta=M\circ e$), which acts on generators of $\mathbb{H}_*(LS^3)$ as follows \cite{T}, \cite{Me}: 
\begin{align*}
\Delta (\alpha \otimes y^i) & =i(1\otimes y^{i-1})\nonumber \\
\Delta (1\otimes y^i) & =0.\nonumber 
\end{align*}
The Gysin exact sequence:\\
\begin{tikzpicture}[descr/.style={fill=white,inner sep=1.5pt}]
        \matrix (m) [
            matrix of math nodes,
            row sep=2em,
            column sep=2.5em,
            text height=1.5ex, text depth=0.25ex
        ]
        { & H^{S^1}_4(LS^3) & H^{S^1}_2(LS^3) & \mathbb{H}_{0}(LS^3)\cong \ZZ (1\otimes 1) \\
          & H^{S^1}_3(LS^3) & H^{S^1}_1(LS^3) & \mathbb{H}_{-1}(LS^3)\cong \ZZ (\alpha \otimes  y) \\
            & H^{S^1}_2(LS^3) & H^{S^1}_0(LS^3) & \mathbb{H}_{-2}(LS^3)\cong 0 \\
            & H^{S^1}_1(LS^3) & 0 & \mathbb{H}_{-3}(LS^3)\cong \ZZ (\alpha \otimes 1) \\
            & H^{S^1}_0(LS^3) &0  & 0 \\
        };

        \path[overlay,->, font=\scriptsize,>=latex]

        (m-1-2) edge (m-1-3)
        (m-1-3) edge[orange] node[above] {$M$}(m-1-4) 
        (m-1-4) edge[out=355,in=175,cyan] node[descr,yshift=0.3ex] {$e$} (m-2-2)
        (m-2-2) edge (m-2-3)
        (m-2-3) edge[orange] node[above] {$M$}(m-2-4)
        (m-2-4) edge[out=355,in=175,cyan] node[descr,yshift=0.3ex] {$e$} (m-3-2)
        (m-3-2) edge (m-3-3)
        (m-3-3) edge[orange] node[above] {$M$} (m-3-4)
        (m-3-4) edge[out=355,in=175,cyan] node[descr,yshift=0.3ex] {$e$} (m-4-2)
        (m-4-2) edge (m-4-3)
        (m-4-3) edge[orange] node[above] {$M$} (m-4-4)
         (m-4-4) edge[out=355,in=175,cyan] node[descr,yshift=0.3ex] {$e$} (m-5-2)   
        (m-5-2) edge (m-5-3)
        (m-5-3) edge[orange] node[above] {$M$}(m-5-4);
\end{tikzpicture}
\\ \indent The short exact sequence in the last two rows shows that $H^{S^1}_0(LS^3)\cong \ZZ (\alpha \otimes 1)$. The short exact sequence in the third and fourth row,
\begin{equation}
0\rightarrow H_1^{S^1}(LS^3) \rightarrow 0\nonumber 
\end{equation} shows that $H^{S^1}_1(LS^3)\cong 0$. Thus, we obtain a short exact sequence from the second and third row,
\begin{equation}
0\rightarrow \ZZ (\alpha \otimes y) \rightarrow H^{S^1}_2(LS^3) \rightarrow \ZZ (\alpha \otimes 1) \rightarrow 0. \nonumber 
\end{equation}
Since the last non-zero term in the sequence is free, the sequence splits, giving $H^{S^1}_2(LS^3)\cong \ZZ (\alpha \otimes y) \oplus \ZZ (\alpha \otimes 1)$. To calculate $H^{S^1}_3(LS^3)$, we use the BV operator. The injective map $e$ in the exact sequence (1.2) means that $e(\alpha \otimes y)= \alpha \otimes y$. Since $\Delta(\alpha \otimes y)=M \circ e (\alpha \otimes y)=1 \otimes 1$, the map $M$ in the first row of the Gysin sequence above is surjective, so the connecting map $e$ from the first to the second row of the Gysin sequence has kernel $\ZZ (1 \otimes 1) $. Thus we have a short exact sequence,
\begin{equation}
0  \rightarrow H^{S^1}_3(LS^3) \rightarrow 0\nonumber 
\end{equation}
So $H^{S^1}_3(LS^3) \cong 0$.\\
\indent We may be able to continue computing the rest of the equivariant homology groups of $LS^3$ in this way, but we eventually reach extension issues. In summary, we have the following remark:
\rem
\begin{eqnarray}
& H^{S^1}_0(LS^3) &=\ZZ (\alpha \otimes 1)\nonumber  \\
& H^{S^1}_1(LS^3) &=0\nonumber \\
& H^{S^1}_2(LS^3) &=\ZZ (\alpha \otimes y) \oplus \ZZ (\alpha \otimes 1)\nonumber  \\
& H^{S^1}_3(LS^3) &=0.\nonumber 
\end{eqnarray}
\normalfont
\subsection{Equivariant Cohomology of $LS^3$}
~
\indent Consider the fibration
\begin{equation}
LS^3 \rightarrow  LS^3\times_{S^1}ES^1 \rightarrow \CC P^{\infty}\nonumber 
\end{equation} and the cohomology Leray-Serre spectral sequence associated with it.\\
\indent We use the fact that we know the ordinary cohomology of $LS^3$ and $\CC P^{\infty}$, since the $E_{\infty}$ page converges to $H^*_{S^1}(LS^3)$, the equivariant cohomology of $S^3$.
\rem \cite{CJY}
\begin{equation}
H^*(LS^3)\cong H^*(\Omega S^3)\otimes H^*(S^3)\cong \Gamma[y]\otimes \Lambda [a] \hspace*{.5cm}|a|=3, |y|=2, y_i=\frac{y^i}{i!} \nonumber 
\end{equation}
\rem
\begin{equation}
H^*(\CC P^{\infty})\cong \ZZ [x], \;|x|=2\nonumber 
\end{equation} 
\normalfont
Below is the $E_2$ page of the spectral sequence. All of the nonzero entries are $\ZZ$ generated by the entry. The arrows are the $d_2$ maps.\\
{\small 
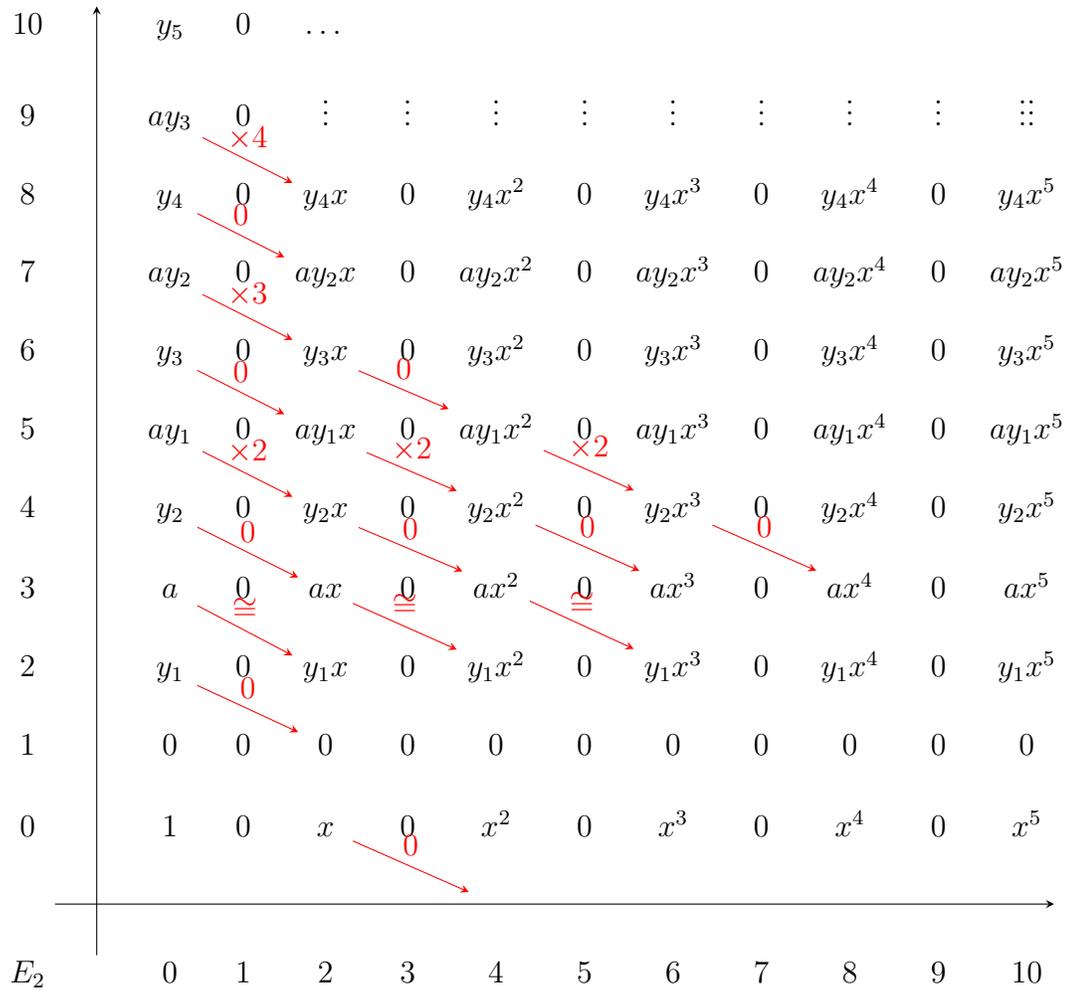
\begin{figure}
\begin{tikzpicture}
  \matrix (m) [matrix of math nodes,
    nodes in empty cells,nodes={minimum width=4ex,
    minimum height=4ex,outer sep=0pt},
    column sep=1ex,row sep=1ex]{
    10&& y_5  &0&\ldots&&&&&&&&&&\\
    9 && ay_3 &0&\vdots&\vdots&\vdots&\vdots&\vdots&\vdots&\vdots&\vdots&\vdots\vdots&&&&&&&&&&\\
    8 && y_4 &0&y_4x&0&y_4x^2&0&y_4x^3&0&y_4x^4&0&y_4x^5&&\\
    7 && ay_2 &0&ay_2x&0&ay_2x^2&0&ay_2x^3&0&ay_2x^4&0&ay_2x^5&&\\
    6 && y_3 &0&y_3x&0&y_3x^2&0&y_3x^3&0&y_3x^4&0&y_3x^5&&\\
    5 && ay_1&0&ay_1x&0&ay_1x^2&0&ay_1x^3&0&ay_1x^4&0&ay_1x^5&&\\
    4 && y_2 &0 & y_2x&  0& y_2x^2 & 0 & y_2x^3 & 0&y_2x^4&0&y_2x^5&&\\
    3 && a  & 0 & ax & 0 & ax^2 & 0 & ax^3 & 0 & ax^4 & 0 & ax^5 &&\\
    2 && y_1 & 0 & y_1x & 0 & y_1x^2 & 0 & y_1x^3 & 0 & y_1x^4 & 0 & y_1x^5 &&\\
    1 && 0 &  0 & 0 & 0 & 0 & 0 & 0 & 0 & 0 & 0 & 0 &&\\
    0 && 1 &  0  & x &  0  & x^2 & 0 & x^3 & 0 & x^4 & 0 & x^5 &&\\
      &&   &      &     &     & &&&&&&&&\\
     \strut  E_2  & &  0  &  1  &  2  &  3 & 4 & 5 & 6 & 7 & 8 & 9 & 10 & \strut \\
    };
  \draw[-stealth] (m-13-2.north) -- (m-1-2.north);
  \draw[-stealth] (m-12-1.east) -- (m-12-13.east);
   \draw[-stealth, red] (m-9-3) -- node[above] {$0$} (m-10-5);
   \draw[-stealth,red] (m-8-3) -- node[above] {$\cong$}(m-9-5);
   \draw[-stealth,red] (m-7-3) -- node[above] {$0$}(m-8-5);
   \draw[-stealth,red] (m-6-3) -- node[above] {$\times 2$}(m-7-5);
   \draw[-stealth,red] (m-5-3) -- node[above] {$0$}(m-6-5);
   \draw[-stealth,red] (m-4-3) -- node[above] {$\times 3$}(m-5-5);
   \draw[-stealth,red] (m-3-3) -- node[above] {$0$}(m-4-5);
   \draw[-stealth,red] (m-2-3) -- node[above] {$\times 4$}(m-3-5);
   \draw[-stealth,red] (m-11-5) --  node[above] {$0$}(m-12-7);
   \draw[-stealth,red] (m-5-5) -- node[above] {$0$}(m-6-7);
   \draw[-stealth,red] (m-6-7) -- node[above] {$\times 2$}(m-7-9);
   \draw[-stealth,red] (m-7-9) -- node[above] {$0$}(m-8-11);
    \draw[-stealth,red] (m-6-5) -- node[above] {$\times 2$}(m-7-7);
    \draw[-stealth,red] (m-7-7) -- node[above] {$0$}(m-8-9);
       \draw[-stealth,red] (m-7-5) -- node[above] {$0$}(m-8-7);
          \draw[-stealth,red] (m-8-7) -- node[above] {$\cong$}(m-9-9);
          \draw[-stealth,red] (m-8-5) -- node[above] {$\cong$}(m-9-7);
\end{tikzpicture}
\caption{The $E_2$ page of the spectral sequence.}
\end{figure}}
\indent We can figure out the first few equivariant cohomology groups easily. It can immediately be seen that $H_{S^1}^0(LS^3)\cong \ZZ (1\otimes 1)$ and $H_{S^1}^1(LS^3)\cong 0$. For $H_{S^1}^2(LS^3)$, the differential maps $d_2$ going to and from the generators along the diagonal, $y_1$ and $x$, are zero, so these generators survive to the $E_{\infty}$ page. In the filtration of $H_{S^1}^2(LS^3)$ corresponding to this spectral sequence, we obtain $0 \subset \ZZ y_1 \subset H_{S^1}^2(LS^3)$
where $H_{S^1}^2(LS^3)/\ZZ x \cong \ZZ y_1$. Thus $H_{S^1}^2(LS^3) \cong \ZZ y_1 \oplus \ZZ x$.\\
\indent The derivation property of the differentials in the Serre spectral sequence makes the computation of the $d_2$ differentials easier. We only need to know the image of $x$, $y_1$, and $a$ through $d_2$ to know the image of the other generators in the $E^2$ grid. We see immediately that $d_2 (x)=0$ and $d_2 (y_1) = 0$. From the multiplicative property of the sequence, we can conclude that $d_2(x^i) =0$ and $d_2(y_i) =0$ for all $i\geq 1$. We computed $H^{S^1}_3(LS^3)\cong 0 $ above, and using the Universal Coefficient Theorem, we find that $H_{S^1}^3(LS^3)\cong 0 $ also. This means that on the $E_{\infty}$ page of the spectral sequence, there should only be zeros along the third diagonal. This gives that $d_2: \ZZ a \rightarrow \ZZ y_1x$ should be an isomorphism. Since $a$ and $y_1x$ are the generators of these isomorphic groups $d_2(a)=\pm y_1x$. Let us assume $d_2(a)=y_1x$. Also,
\begin{equation}
d_2(ay_i)=d(a)y_i=y_1xy_i=y_1x\frac{y_1^i}{i!}=(i+1)y_{i+1}x\nonumber 
\end{equation} To summarize:
\begin{eqnarray}
& d_2(x^i) &=0\nonumber \\
& d_(y_i) &=0\nonumber \\
& d_2(a)&=y_1x\nonumber \\
& d_2(ay_i) &=(i+1)y_{i+1}x.\nonumber 
\end{eqnarray}\\
These calculations correspond to the red arrows on the $E_2$ page above. The spectral sequence collapses at the $E_3$ page since there can never be nonzero differentials after the $E_2$ page because there is nothing for these differentials to hit, so $E_3=E_{\infty}$.\\
\indent Let's take a look at the $E_{\infty}$ page:\\

\begin{figure}{\footnotesize
 \begin{tikzpicture}
  \matrix (m) [matrix of math nodes,
    nodes in empty cells,nodes={minimum width=4ex,
    minimum height=4ex,outer sep=0pt},
    column sep=1ex,row sep=1ex]{
    9 &&\vdots& \vdots&\vdots&\vdots&\vdots&\vdots&\vdots&\vdots&\vdots&\vdots&&&&\\
    8 && y_4 &0&\ZZ y_4x/4\ZZ&0&\ZZ y_4x^2/4\ZZ&0&\ZZ y_4x^3/4\ZZ&0&\ZZ y_4x^4/4\ZZ&\ldots&&\\
    7 && 0&0&0&0&0&0&0&0&0&\ldots&&\\
    6 && y_3 &0&\ZZ y_3x/ 3\ZZ&0&\ZZ y_3x^2/ 3\ZZ&0&\ZZ y_3x^3\ 3\ZZ&0&\ZZ y_3x^4/3\ZZ&\ldots&&\\
    5 && 0&0&0&0&0&0&0&0&0&\ldots&&\\
    4 && y_2 &0 & \ZZ y_2x/2\ZZ&  0& \ZZ y_2x^2/ 2\ZZ & 0 & \ZZ y_2x^3/ 2\ZZ & 0& \ZZ y_2x^4/ 2\ZZ&\ldots&&\\
    3 && 0 & 0 & 0 & 0 & 0 & 0 & 0& 0 & 0 & \ldots &&\\
    2 && y_1 & 0 &0& 0 & 0& 0 & 0 & 0 & 0 & \ldots&&\\
    1 && 0 &  0 & 0 & 0 & 0 & 0 & 0 & 0 & 0 & \ldots &&\\
    0 && 1 &  0  & x &  0  & x^2 & 0 & x^3 & 0 & x^4 & \ldots &&\\
      &&   &      &     &     & &&&&&&&&\\
     \strut  E_{\infty}  & &  0  &  1  &  2  &  3 & 4 & 5 & 6 & 7 & 8 &  & \strut \\
    };
  \draw[-stealth] (m-12-2.north) -- (m-1-2.north);
  \draw[-stealth] (m-11-1.east) -- (m-11-12.east);
  \end{tikzpicture}}
  \caption{The $E_\infty$ page of the spectral sequence.}
  \end{figure}
  \newpage
  \indent After some work, and with the assumption that $H_{S^1}^i(LS^3)$ is just the direct sum of the diagonal on the $E_\infty$ page shown above, we can make the following remark.
  \rem 
  \begin{align*}
   H_{S^1}^{2i+1}(LS^3) &=0\nonumber \\
    H_{S^1}^{2i}(LS^3) &=\ZZ y_i \oplus \ZZ x^i \oplus \sum_{j=2}^{i-1} \ZZ y_j x^{i-j} /j\ZZ, \hspace*{.5cm} i>0\nonumber \\
    H^0_{S^1}(LS^3) &=\ZZ 1\nonumber .
  \end{align*}

\normalfont

\subsection{Equivariant Homology of $LS^3$}
~
\indent Using the results above and the Universal Coefficient Theorem, we get the following:

\thm
\begin{align}
H_{2i}^{S^1}(LS^3)&=\ZZ(\alpha \otimes y^i) \oplus \ZZ x_i ,\hspace*{.5cm} i \geq 0\nonumber \\
H_{2i+1}^{S^1}(LS^3)&=\sum_{j=2}^{i}\ZZ(\alpha \otimes y^i)x/j\ZZ , \hspace*{.5cm} i \geq 2\nonumber \\
H_1^{S^1}(LS^3)&=H_3^{S^1}(LS^3)=0.\nonumber 
\end{align}
\normalfont
Note that $\alpha \otimes y^i$ is dual to $y_i$ and $1\otimes y^i$ is dual to $ay_i$. This matches the findings of Basu \cite{Ba} and Westerland \cite{We} using the Universal Coefficient Theorem.



\subsection{The Spectral Sequence Associated with the Gysin Sequence}
~
\indent To determine the erasing and marking maps, we will translate the Gysin sequence into a spectral sequence and see how they arise in the computation of the spectral sequence. We are using the fact that $H^*(S^1)=\Lambda \omega$ where $|\omega|=1$. Note that on the $E_{\infty}$ page, $a=\omega y_1$.\\

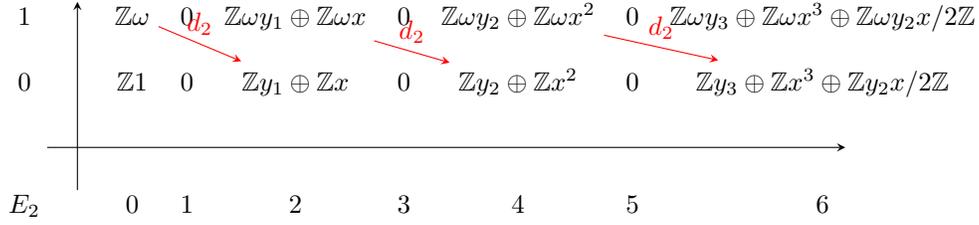
\begin{figure}
{\footnotesize
\begin{tikzpicture}
  \matrix (m) [matrix of math nodes,
    nodes in empty cells,nodes={minimum width=4ex,
    minimum height=4ex,outer sep=0pt},
    column sep=.5ex,row sep=1ex]  
{    1 && \ZZ \omega &  0 & \ZZ \omega y_1 \oplus \ZZ \omega x & 0 & \ZZ \omega y_2 \oplus \ZZ \omega x^2 & 0 & \ZZ \omega y_3 \oplus \ZZ \omega x^3 \oplus \ZZ \omega y_2 x/2\ZZ & \\
  0 && \ZZ 1 &  0 & \ZZ  y_1 \oplus \ZZ  x & 0 & \ZZ  y_2 \oplus \ZZ  x^2 & 0 & \ZZ  y_3 \oplus \ZZ  x^3 \oplus \ZZ y_2 x/2\ZZ &\\
      &&   &      &     &     & &&&&&&&\\
     \strut  E_{2}  & &  0  &  1  &  2  &  3 & 4 & 5 & 6 &  \strut \\
    };
  \draw[-stealth] (m-4-2.north) -- (m-1-2.north);
  \draw[-stealth] (m-3-1.east) -- (m-3-9.east);
  \draw[-stealth, red] (m-1-3) -- node[above] {$d_2$} (m-2-5);
  \draw[-stealth, red] (m-1-5) -- node[above] {$d_2$} (m-2-7);
  \draw[-stealth, red] (m-1-7) -- node[above] {$d_2$} (m-2-9);;
  
  \end{tikzpicture}}
 \caption{The $E_2$ page.} 
\end{figure}  
  \vspace*{2cm}
  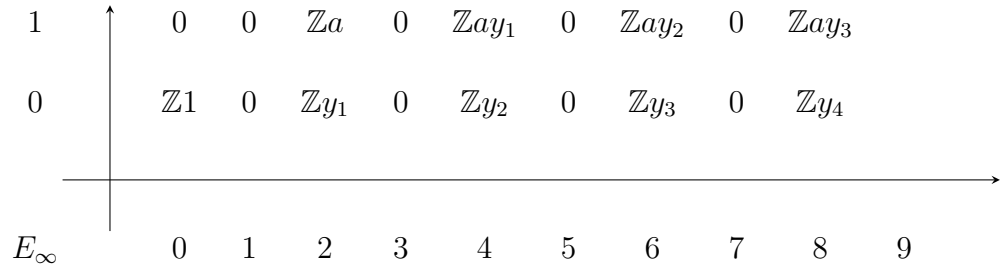
\begin{figure}
  \begin{tikzpicture}
  \matrix (m) [matrix of math nodes,
    nodes in empty cells,nodes={minimum width=4ex,
    minimum height=4ex,outer sep=0pt},
    column sep=1ex,row sep=1ex]  
{    1 && 0 &  0 & \ZZ a & 0 & \ZZ ay_1 & 0 & \ZZ ay_2 & 0 & \ZZ ay_3 &&\\
    0 &&\ZZ 1 &  0 & \ZZ y_1 & 0 & \ZZ y_2 & 0 & \ZZ y_3 & 0 & \ZZ y_4 && \\
      &&   &      &     &     & &&&&&&&\\
     \strut  E_{\infty}  & &  0  &  1  &  2  &  3 & 4 & 5 & 6 & 7 & 8 & 9 &  \strut \\
    };
  \draw[-stealth] (m-4-2.north) -- (m-1-2.north);
  \draw[-stealth] (m-3-1.east) -- (m-3-13.east);

  \end{tikzpicture}
  \caption{The $E_\infty$ page.}
  \end{figure}
   \indent First, we will compute the differential maps for the fibration (1.1) in the spectral sequence. Note that the $E_{\infty}=\cdots =E_4=E_3$ because all of the differential maps $d_i$ for $i\geq 3$ are $0$. \\
   \indent Since $H^1(LS^3)=0$, the $E_{\infty}$ page has zeros along the $1$ diagonal. This means for the map $d_2: \ZZ \omega \rightarrow \ZZ y_1 \oplus \ZZ x$, $ker(d_2)=0$ so the map is injective. Since the entry $E_{\infty}^{2,0}=\ZZ y_1$, the image of $d_2$ must be $\ZZ x$, so we can say that $d_2(\omega)=x$ (up to a sign). We have that $d_2(y_1)=0$ and $d_2(x)=0$ also. Using the multiplicative structure of the spectral sequence, we obtain the following remark. \\

\rem
\begin{align}
d_2(\omega y_i) & =xy_i\nonumber \\
d_2(\omega x^i)&=x^{i+1}\nonumber \\
d_2(\omega y_i x^j) & = y_i x^{j+1}\nonumber 
\end{align}
\normalfont


\subsection{The Erasing Map $e$}
~
\indent In the Gysin sequence above for the circle fibration~\ref{eqn:threecirclebundle}, the erasing map $e:$ can be viewed as the map induced  by $\epsilon$, the projection map in the fibration~\ref{eqn:threecirclebundle}, so $e=\epsilon_*$. Since we are interested in seeing how the erasing map acts on specific generators, we will instead look at the dual map $\epsilon^*$.\\
\indent The map $\epsilon^*:H_{S^1}^i(LS^3) \rightarrow H^i(LS^3)$ can be derived from the spectral sequence of the above fibration~\ref{eqn:threecirclebundle}. The map $\epsilon^*$ is the composition of the surjection map $E_2^{i,0}=H_{S^1}^i(LS^3) \rightarrow E_{\infty}^{i,0}=E_2^{i,0}/im(d_2)$ and the inclusion map $E_{\infty}^{i,0}\rightarrow H^i(LS^3)$ \cite{Mc}. This means that an image of a generator in the bottom row of the $E^2$ page of the spectral sequence by $\epsilon^*$ is nonzero if it survives to the $E^{\infty}$ page, and a generator's image is zero if it does not survive. The following remark is immediate.
\rem 
\begin{align}
\epsilon ^*(x^i)&=0 \nonumber \\
\epsilon ^*(y_i)&=y_i \nonumber \\
\epsilon ^*(y_jx^i)&=0\nonumber\\
\end{align}\\
\normalfont
To dualize $\epsilon^*$ to obtain $e$, we need the Kronecker pairing as in the computation of $M$.
\lem
\begin{align}
e=\epsilon _* :H_*(LS^3) & \rightarrow H_*^{S^1}(LS^3)\nonumber \\
\alpha \otimes y^i & \mapsto \alpha \otimes y^i\nonumber \\
 1\otimes y^i & \mapsto  (\alpha \otimes y^{i+1})x \hspace*{.25cm}\nonumber \\
1 \otimes 1 & \mapsto 0. \nonumber
\end{align}
\proof Since $|\alpha \otimes y^i |=2i$, $e(\alpha \otimes y^i )=kx_i+l(\alpha \otimes y^i) $ for $k,l\in \ZZ$. We have,
\begin{align*}
<\epsilon ^*(y_i),\alpha \otimes y^i >&=<y_i,\alpha \otimes y^i >=1 
\\& =<y_i, \epsilon_*(\alpha \otimes y^i )> \\
&=<y_i,kx_i+l(\alpha \otimes y^i )>\\
&=k<y_i,x_i>+l<y_i,\alpha \otimes y^i >=l 
\end{align*} and 
\begin{align*}
<\epsilon ^*(x^i),\alpha \otimes y^i >&=<0,\alpha \otimes y^i >=0\\
&=<x^i, \epsilon_*(\alpha \otimes y^i )>\\
 &=<y^i,kx_i+\alpha \otimes y^i >\\
 &=k<x^i,x_i>+l<x^i,\alpha \otimes y^i >=k. \nonumber
\end{align*}
Therefore $e(\alpha \otimes y^i )=\epsilon_*(\alpha \otimes y^i )=\alpha \otimes y^i $.\\
\indent Since $H_3^{S^1}(LS^3)=0$ and $|1\otimes 1|=3$, we must have $e(1\otimes 1)=0$. Since $|1 \otimes y^i |=3+2i$, these generators are of odd degree so they cannot be paired with generators in cohomology since $H^{2i+1}_{S^1}(LS^3) =0$ for $i \geq 0$, so we need to use another technique to find the image of $1 \otimes y^i$. For this we will go back to the Gysin sequence for the fibration (1.1).\\
\indent For $i=1$ we look at the following piece of the Gysin sequence.\\
\begin{tikzpicture}[descr/.style={fill=white,inner sep=1.5pt}]
        \matrix (m) [
            matrix of math nodes,
            row sep=2em,
            column sep=2em,
            text height=1.5ex, text depth=0.25ex
        ]
        { & H^{S^1}_4(LS^3))\cong \ZZ(\alpha \otimes y^2) \oplus \ZZ x_2 & \mathbb{H}_{2}(LS^3)\cong \ZZ (1\otimes y) &\\
         & H^{S^1}_5(LS^3) \cong \ZZ /2 \ZZ ((\alpha \otimes y^2)x) & 0 \\
             };

        \path[overlay,->, font=\scriptsize,>=latex]

        (m-1-2) edge[orange] node[above] {$M$}(m-1-3) 
        (m-1-3) edge[out=355,in=175,cyan] node[descr,yshift=0.3ex] {$e$} (m-2-2)
        (m-2-2) edge (m-2-3)
       ;
\end{tikzpicture}\\
Since $e$ is surjective and $ker(e)=im(M)\cong 2\ZZ(1\otimes y)$ by Lemma 6.2, we have that $im(e)\cong \ZZ/ 2\ZZ (1\otimes y)$. Thus $e(1\otimes y)=\pm(\alpha \otimes y^2)x$.\\
\indent In general, we have the following.
\vspace*{.5cm}\\
\begin{tikzpicture}[descr/.style={fill=white,inner sep=.5pt}]
        \matrix (m) [
            matrix of math nodes,
            row sep=2em,
            column sep=.5em,
            text height=1.5ex, text depth=0.25ex
        ]
        { &  H^{S^1}_{2i+2}(LS^3)\cong \ZZ(\alpha \otimes y^{i+1}) \oplus \ZZ x_{i+1} & \mathbb{H}_{2i}(LS^3)\cong \ZZ (1\otimes y^i) & \\
          & ~~~~H^{S^1}_{2i+3}(LS^3)\cong \sum_{j=2}^{i+1}\ZZ_j(\alpha\otimes y^j)x_{i-j
        +2} 
          &  & 
          \\ &  H^{S^1}_{2i+1}(LS^3) \cong \sum_{j=2}^{i}\ZZ_j(\alpha\otimes y^j)x_{i-j
        +1} & 0\\
             };

        \path[overlay,->, font=\scriptsize,>=latex]

        (m-1-2) edge[orange] node[above] {$M$}(m-1-3) 
        (m-1-3) edge[out=354,in=174,cyan] node[descr,yshift=0.3ex] {$e$} (m-2-2)
       
        (m-2-2) edge[orange, out=355,in=175] node[descr,yshift=0.3ex] {$c$}  (m-3-2)  
         (m-3-2)edge[black] node[above] {$M$}(m-3-3) 
       ;
\end{tikzpicture}\\
We have that $ker(e)=im(M)=(i+1)\ZZ (1\otimes y^i)$, so $im(e)\cong (\ZZ_{i+1} \ZZ )(1\otimes y^i)$. Since the map $c$ is given by the cap product with $x\in  H^2_{S^1}(LS^3)$, $c((\alpha \otimes y^j)x_{i-j+2})=(\alpha \otimes y^j)x_{i-j+1}$. So $ker(c)=im(e) \cong \ZZ_j (\alpha \otimes y^{i+1})x$.  We can conclude that $e(1 \otimes y^i)=(\alpha \otimes y^{i+1})x$.
\qed

\subsection{The Marking Map $M$}
~
\indent We consider the dual of the marking map, $M^*:H^i(LS^3) \rightarrow H^{i-1}_{S^1}(LS^3)$. This can be derived from the spectral sequence of the circle fibration~\ref{eqn:threecirclebundle}. $M_*$ is the composition of the surjective map $H^i(LS^3)\longrightarrow E^{i-1,1}_{\infty}\cong H^i(LS^3)/E_{\infty}^{i,0}$ and the injective map $E^{i-1,1}_{\infty}\cong ker(d_2) \longrightarrow E^{i-1,1}_{2}$ \cite{Mc}.
\rem
\begin{align*}
M^*(y_i)&=0\nonumber \\
M^*(ay_i) &= (i+1)y_{i+1}\nonumber 
\end{align*}
\begin{proof} Since the kernel of the differential $d_2: (E_2^{-1,1}\cong 0) \longrightarrow (E_2^{1,0}\cong 0)$ is $0$, for $1\in H^0(LS^3)$, $M^*(1)=0$. To find the image of $y_i$, we consider the composition $M^*:(H^{2i}(LS^3)\cong \ZZ y_i) \longrightarrow (E_{\infty}^{2i-1,1}\cong H^{2i}(LS^3)/E_{\infty}^{2i,0} \cong \ZZ y_i/\ZZ y_i\cong 0)$, so $M^*(y_i)=0$. The image of $ay_i$ can be determined by identifying $ay_i$ with $(i+1)\omega y_{i+1}$. Then,
\begin{align}
M^*: H^{2i+1}  &\rightarrow  E_{\infty}^{2i+2,1} \cong ker(d_2)  \longrightarrow   E_2^{2i+2,1} \nonumber \\
ay_i  &\mapsto  (i+1)\omega y_{i+1}  \xmapsto{\phantom{ E_2^{2i+2,1} }} (i+1)\omega y_{i+1}=(i+1) y_{i+1}\nonumber 
\end{align}
so $M^*(ay_i)=(i+1)y_{i+1}$.
\end{proof}
\normalfont
To dualize $M_*$ to obtain $M$, we need the Kronecker pairing \cite{Br}.

\defn The Kronecker pairing is a map
\begin{align*}
<-,->:H^i(X)\otimes H_i(X) \rightarrow \ZZ \nonumber 
\end{align*}
such that for $\alpha=[f] \in H^i(X)$ and $\gamma=[c]\in H_i(X)$ then 
\begin{align*}
<\alpha , \gamma >=f(c).\nonumber 
\end{align*}
Alternatively, for $\beta :H^i(X) \rightarrow Hom(H_i(X))$, the map from the universal coefficient theorem, 
\begin{equation}
<\alpha , \gamma >=f(c) \in \ZZ . \nonumber 
\end{equation}\\
\normalfont
The Kronecker pairing satisfies the following property, which will be used to dualize the map $M ^*$:
\begin{equation}
<f^*(\alpha),\gamma >=<\alpha , f_*(\gamma )>\nonumber 
\end{equation}
\lem
\begin{align*}
M=M_*: H_*^{S^1}(LS^3)& \longrightarrow  H_{*+1}(LS^3)\nonumber \\
\alpha \otimes 1  & \xmapsto{\phantom{ H_*^{S^1}(LS^3)}}  0\nonumber \\
\alpha \otimes y^i & \xmapsto{\phantom{ H_*^{S^1}(LS^3)}}  i(1\otimes y^{i-1})\nonumber \\
x_i & \xmapsto{\phantom{ H_*^{S^1}(LS^3)}} 0\nonumber \\
(\alpha \otimes y^j)x_i & \xmapsto{\phantom{ H_*^{S^1}(LS^3)}}  0\nonumber 
\end{align*}
\proof Since $H_1(LS^3)=0$, it is immediate that $M(\alpha \otimes 1)=0$. To find the image of $\alpha \otimes y^i\in H_{2i}^{S^1}(LS^3)$, since $M$ is a map of degree $+1$, the only possible generator of $H_*(LS^3)$ of degree $2i+1$ is $ 1\otimes y^{i-1}$. Let $M_*(\alpha \otimes y^i)=k (1\otimes y^{i-1})$ for $k\in \ZZ$. Then
\begin{align*}
<M^*( ay_{i-1}), y^i>&=<ay_{i-1},M_*(\alpha \otimes y^i)>\\
&=<ay_{i-1},k(1 \otimes y^{i-1})>=k\\
& =<iy_i,\alpha \otimes y^i>=i.
\end{align*}
This implies $k=i$, so $M(\alpha \otimes y^i)=i(1\otimes y^{i-1})$.\\
\indent We must have $M(x_i)=0$ since $|x_i|=2i$ and the only generator on cohomology that it can be paired with is $ay^{i-1}$, which is not dual to $x_i$. Similarly, $(\alpha \otimes y^j)x_i$ gets sent to zero by $M$ since it is torsion, mapping into a free group. 
\qed
\normalfont

\subsection{The String Bracket $[-,-]$}
~
\indent Recall the string bracket from Definition~\ref{defn:stringbracket}
\begin{align*}
[a,b]=(-1)^{(|a|-3})e(M(a)\bullet M(b))
\end{align*}
of degree $-1$ for $LS^3$.\\
\indent The only possible non-zero bracket is from the pair $\alpha \otimes y^i, \alpha \otimes y^j$, as M maps all other generators of $H_*^{S^1}(LS^3)$ to $0$, thus the string bracket of these generators are also $0$. We see that when $i\geq 1$ and $j \geq 1$,
\begin{align*}
\nonumber [\alpha \otimes y^i,\alpha \otimes y^j]&=(-1)^{2i-3}e(M(\alpha \otimes y^i)\bullet M( \alpha \otimes y^j))\\
&=-e(i(1\otimes y^{i-1}) \cdot j(1 \otimes y^{j-1})) \\
&=-e((ij)(1\otimes y^{i+j-2}))\\
&=-ij(\alpha \otimes y^{i+j-1})x. 
\end{align*} 
So the bracket is equal to zero if both $i=1$ and $j=1$ or if $(i+j-1)\mid ij$ and non-zero in all other cases.\\
\indent As it turns out, the only non-zero brackets are torsion elements, which corresponds to the findings of \cite{Ba}, which are that the brackets are all trivial when considering rational string homology of $S^3$.


\section{String Homology and the String Bracket of Odd Spheres}

\subsection{String Homology for Odd Spheres}
~
\indent We try to compute the string homology for odd spheres using only the Gysin sequence for the following fibration:
\begin{align} \label{eqn:oddcirclebundle} 
 S^1 \rightarrow LS^n \times ES^1 \rightarrow LS^n \times_{S^1} ES^1 
\end{align}
for $n$ odd.\\
\indent Recall that the loop homology is given as follows by \cite{CJY}:
\begin{equation}
\mathbb{H}_*(LS^n)=\Lambda[a]\otimes \ZZ [u] \nonumber
\end{equation}
where $a$ corresponds to the dual of $[S^n]$, so $|a|=-n$ and $|u|=n-1$ after re-grading. \\
\indent The BV-operator acts on the generators as follows, \cite{Me}:
\begin{align*}
\Delta(a\otimes u^i)= & i(1 \otimes u^{i-1})\\
\Delta (1 \otimes u^i)= & 0.
\end{align*}
Throughout this section, we consider $n$ to be odd.
\normalfont
Consider the bottom of the Gysin sequence. Let $H^{S^1}_i$ denote $H^{S^1}_i(LS^n)$ and $\mathbb{H}_i$ denote $\mathbb{H}_i(LS^n)$\\
\begin{tikzpicture}[descr/.style={fill=white,inner sep=1.5pt}]
        \matrix (m) [
            matrix of math nodes,
            row sep=2em,
            column sep=2.5em,
            text height=1.5ex, text depth=0.25ex
        ]
        { & H^{S^1}_{n+3} & H^{S^1}_{n+1} & \mathbb{H}_{2}\cong 0 \\        
        & H^{S^1}_{n+2} & H^{S^1}_n\cong 0 & \mathbb{H}_{1}\cong 0\\                  
             & H^{S^1}_{n+1} & H^{S^1}_{n-1} \cong \ZZ \oplus \ZZ& \mathbb{H}_{0}\cong \ZZ (1\otimes 1) \\     
        & H^{S^1}_{n}\cong 0 & H^{S^1}_{n-2}\cong 0 & \mathbb{H}_{-1}\cong \ZZ (a\otimes u) \\     
        & H^{S^1}_{n-1} & H^{S^1}_{n-3}\cong \ZZ (\gamma_{\frac{n-3}{2}}) & \mathbb{H}_{-2}\cong 0 \\     
        & H^{S^1}_{n-2}\cong 0 & H^{S^1}_{n-4} & 0\\        
        &\vdots &\vdots &\vdots\\
        & H^{S^1}_4 \cong \ZZ (\gamma_2) & H^{S^1}_2\cong \ZZ (\gamma) & \mathbb{H}_{-n+3}\cong 0 \\
        & H^{S^1}_3 & H^{S^1}_1 & \mathbb{H}_{-n+2}\cong 0 \\
          & H^{S^1}_2\cong \ZZ( \gamma) & H^{S^1}_0\cong \ZZ & \mathbb{H}_{-n+1}\cong 0 \\
          & H^{S^1}_1\cong 0 & 0 & \mathbb{H}_{-n}\cong \ZZ (a \otimes 1) \\
            & H^{S^1}_0\cong \ZZ &0  &  \\
        };

        \path[overlay,->, font=\scriptsize,>=latex]

        (m-1-2) edge  node[above] {$\cong$}(m-1-3)
        (m-1-3) edge[orange] node[above] {$M$}(m-1-4) 
        (m-1-4) edge[out=355,in=175,cyan] node[descr,yshift=0.3ex] {$e$} (m-2-2)
        (m-2-2) edge node[above] {$\cong$} (m-2-3)
        (m-2-3) edge[orange] node[above] {$M$}(m-2-4)
        (m-2-4) edge[out=355,in=175,cyan] node[descr,yshift=0.3ex] {$e$} (m-3-2)
        (m-3-2) edge (m-3-3)
        (m-3-3) edge[orange] node[above] {$M$} (m-3-4)
        (m-3-4) edge[out=355,in=175,cyan] node[descr,yshift=0.3ex] {$e$} (m-4-2)
        (m-4-2) edge (m-4-3)
        (m-4-3) edge[orange] node[above] {$M$} (m-4-4)
         (m-4-4) edge[out=355,in=175,cyan] node[descr,yshift=0.3ex] {$e$}(m-5-2)   
        (m-5-2) edge (m-5-3)
        (m-5-3) edge[orange] node[above] {$M$}(m-5-4)
         (m-5-4) edge[out=355,in=175,cyan] node[descr,yshift=0.3ex] {$e$}(m-6-2)   
        (m-6-2) edge  node[above] {$\cong$} (m-6-3)
        (m-6-3) edge[orange] node[above] {$M$}(m-6-4)
        (m-8-2) edge node[above] {$\cong$} (m-8-3)
        (m-8-3) edge[orange] node[above] {$M$}(m-8-4)
        (m-8-4) edge[out=355,in=175,cyan] node[descr,yshift=0.3ex] {$e$}(m-9-2)   
        (m-9-2) edge node[above] {$\cong$} (m-9-3)
        (m-9-3) edge[orange] node[above] {$M$}(m-9-4)
        (m-9-4) edge[out=355,in=175,cyan] node[descr,yshift=0.3ex] {$e$}(m-10-2)   
        (m-10-2) edge node[above] {$\cong$} (m-10-3)
        (m-10-3) edge[orange] node[above] {$M$}(m-10-4)
        (m-10-4) edge[out=355,in=175,cyan] node[descr,yshift=0.3ex] {$e$}(m-11-2)   
        (m-11-2) edge node[above] {$\cong$}(m-11-3)
        (m-11-3) edge[orange] node[above] {$M$}(m-11-4)
        (m-11-4) edge[out=355,in=175,cyan] node[descr,yshift=0.3ex] {$e,\cong$}(m-12-2)   
        (m-12-2) edge (m-12-3);
\end{tikzpicture}\\
The maps $H^{S^1}_i\longrightarrow H^{S^1}_{i-2}$ are given by the cap product with the class generator $x\in H^2(\mathbb{C}P^{\infty})$. Since $H^{S^1}_2(LS^n)\cong \ZZ$, we denote the generator by $\gamma$, which is dual to $x$. We use the notation $\gamma_i=\frac{\gamma^i}{i!}$, dual to $x^i$. Since the maps given by the cap product are isomorphisms between where the loop homology is zero, we have that
\begin{align*}
H^{S^1}_{2i+1}(LS^n)=& 0, \hspace*{.5cm} 0\leq i \leq \frac{n-3}{2}\\
H^{S^1}_{2i}(LS^n)=& \ZZ (\gamma_i), \hspace*{.5cm} 1 \leq i \leq \frac{n-3}{2}.
\end{align*}
Note that for even degrees, the generator $\gamma_i$ increases subscript as isomorphisms in the sequence are given by cap product with $x$, dual to the cup product with $x$.\\
\indent To determine $H^{S^1}_{n-1}(LS^n)$, note that we have a short exact sequence,
\begin{equation}
0\longrightarrow \ZZ (a \otimes u) \longrightarrow H^{S^1}_{n-1} \longrightarrow  \ZZ (\gamma_{\frac{n-3}{2}}) \longrightarrow 0\nonumber
\end{equation}
that splits since the last term is free. Thus $H^{S^1}_{n-1}(LS^n)\cong \ZZ e(a \otimes u) \oplus (\gamma_{\frac{n-1}{2}})$. We use the notation of $e(-)$ to denote that the generator comes from the erasing map. Using the BV-operator to determing the marking map $M:H^{S^1}_{n-1} \rightarrow \mathbb{H}_0$, we have that $M(a \otimes u)=1\otimes 1$, so the erasing map $e:\mathbb{H}_0\rightarrow H^{S^1}_n$ is zero, thus $H_n^{S^1}(LS^n)\cong 0$.
\lem \label{lem:marking} $M(\gamma_{\frac{n-1}{2}})=0$, or more generally, the marking map sends generators coming from $H_*(\mathbb{C}P^{\infty})$ to zero. 
\proof 
In the circle bundle~(\ref{eqn:oddcirclebundle})
the marking map is an umkehr map coming from the projection map. Notice that $\mathbb{C}P^{\infty}=BS^1=\{pt\}\times ES^1 \subset LS^n \times_{S^1} ES^1$. Since $\pi^{-1}(\{pt\}\times_{S^1} ES^1)=\{pt\} \times ES^1$, which is contractible, then $M$ maps generators from $\mathbb{C}P^{\infty}$ into a contractible space, thus $M(\gamma_{i})=0$ for any $i$, where $\gamma_i$ denotes a generator coming from the homology of $\mathbb{C}P^{\infty}$.
\qed \\
\normalfont
\indent With the knowledge that $M(\gamma_{\frac{n-1}{2}})=0$, the cap product map $H^{S^1}_{n+1}\rightarrow H^{S^1}_{n-1}$ is injective with image isomorphic to $\ZZ (\gamma_{\frac{n-1}{2}})$ so $H^{S^1}_{n+1}\cong \ZZ (\gamma_{\frac{n+1}{2}})$.\\
\indent Now consider the next piece of the Gysin sequence where loop homology is non-zero.\\

\begin{tikzpicture}[descr/.style={fill=white,inner sep=1.5pt}]
        \matrix (m) [
            matrix of math nodes,
            row sep=2em,
            column sep=2.5em,
            text height=1.5ex, text depth=0.25ex
        ]
        {         
         H^{S^1}_{2n+1} & H^{S^1}_{2n-1} & \mathbb{H}_{n}\cong 0\\                  
              H^{S^1}_{2n} & H^{S^1}_{2n-2} \cong \ZZ \oplus \ZZ& \mathbb{H}_{n-1}\cong \ZZ (1\otimes u) \\     
         H^{S^1}_{2n-1} & H^{S^1}_{2n-3}\cong 0 & \mathbb{H}_{n-2}\cong \ZZ (a\otimes u^2) \\     
         H^{S^1}_{2n-2} & H^{S^1}_{2n-4}\cong \ZZ (\gamma_{n-2}) & \mathbb{H}_{n-3}\cong 0 \\     
         H^{S^1}_{2n-3}\cong 0 & H^{S^1}_{2n-5}\cong 0 & 0\\        
        \vdots &\vdots &\vdots\\
        H^{S^1}_{n+3} \cong \ZZ (\gamma_{n-2}) & H^{S^1}_{n+1}\cong \ZZ (\gamma_{\frac{n+1}{2}}) & \mathbb{H}_{2}\cong 0 \\
         H^{S^1}_{n+2} & H^{S^1}_n & \mathbb{H}_{1}\cong 0 \\
        };

        \path[overlay,->, font=\scriptsize,>=latex]

        (m-1-1) edge  node[above] {$\cong$}(m-1-2)
        (m-1-2) edge[orange] node[above] {$M$}(m-1-3) 
        (m-1-3) edge[out=355,in=175,cyan] node[descr,yshift=0.3ex] {$e$} (m-2-1)
        (m-2-1) edge  (m-2-2)
        (m-2-2) edge[orange] node[above] {$M$}(m-2-3)
        (m-2-3) edge[out=355,in=175,cyan] node[descr,yshift=0.3ex] {$e$} (m-3-1)
        (m-3-1) edge (m-3-2)
        (m-3-2) edge[orange] node[above] {$M$} (m-3-3)
        (m-3-3) edge[out=355,in=175,cyan] node[descr,yshift=0.3ex] {$e$} (m-4-1)
        (m-4-1) edge (m-4-2)
        (m-4-2) edge[orange] node[above] {$M$} (m-4-3)
         (m-4-3) edge[out=355,in=175,cyan] node[descr,yshift=0.3ex] {$e$}(m-5-1)   
        (m-5-1) edge (m-5-2)
        (m-5-2) edge[orange] node[above] {$M$}(m-5-3)

                (m-7-1) edge node[above] {$\cong$} (m-7-2)
        (m-7-2) edge[orange] node[above] {$M$}(m-7-3)
         (m-7-3) edge[out=355,in=175,cyan] node[descr,yshift=0.3ex] {$e$}(m-8-1)   
        
        (m-8-1) edge node[above] {$\cong$} (m-8-2)
        (m-8-2) edge[orange] node[above] {$M$}(m-8-3)

    ;
\end{tikzpicture}\\
In the third and fourth row above, we have a short exact sequence with $H^{S^1}_{2n-2}$ in the center, which splits, so $H^{S^1}_{2n-2}\cong \ZZ (e(a\otimes u^2)) \oplus \ZZ (\gamma_{n-1})$. Mapping $H^{S^1}_{2n-2}$ through $M$, we have $M(e(a\otimes u^2))=2(1\otimes u)$ given by the BV-operator.
\indent Thus $H^{S^1}_{2n-1}\cong \ZZ/2\ZZ (e(1 \otimes u))$. Since the cap product map $H^{S^1}_{2n}\rightarrow H^{S^1}_{2n-2}$ is injective with image $\ZZ (\gamma_{n-1})$, $H^{S^1}_{2n} \cong \ZZ (\gamma_{n})$. Summarizing, we have
\begin{align*}
H^{S^1}_{2i+1}(LS^n)\cong & 0, \hspace*{.5cm} & \frac{n-1}{2}\leq i \leq n-2 \\
H^{S^1}_{2i}(LS^n)\cong & \ZZ (\gamma_i), \hspace*{.5cm} & \frac{n+1}{2} \leq i \leq n-2  \\
H^{S^1}_{n-2}(LS^n)\cong & \ZZ e(a\otimes u^2) \oplus \ZZ (\gamma_{n-1}) \\
H^{S^1}_{n-1}(LS^n)\cong & \ZZ/2\ZZ (e(1 \otimes u)) \\
H^{S^1}_{2n}(LS^n)\cong & \ZZ (\gamma_{n}) \\
H^{S^1}_{2i+1} \cong & \ZZ_2 (1 \otimes u)\gamma_{i-n-1} \hspace*{.5cm}& n-1 \leq i \leq \frac{3n-5}{2} \\
H^{S^1}_{2i} \cong & \ZZ \gamma_i \hspace*{.5cm} & \frac{n}{2} \leq i \leq \frac{3n-5}{2} 
\end{align*}
Now assume the following holds for all $k\in \NN$:

\begin{equation} \label{eqn:oddinduction}
H^{S^1}_{2i} \cong \left\{
\begin{array}{rl}
\ZZ \gamma_i & \text{if } (n-1)\nmid 2i,\\
\ZZ \gamma_i \oplus \ZZ(e(a\otimes u^i)) & \text{if } (n-1)\vert 2i.
\end{array} \right.
\end{equation}
for $1\leq i \leq \frac{(k+1)(n-1)-2}{2}$ and 
\begin{equation}
H_{2i+1}^{S^1}\cong t_k \nonumber
\end{equation}
for $\frac{k(n-1)}{2}\leq i \leq \frac{(k+1)(n-1)-2}{2}$, where $t_k$ is a torsion group of order $k!$. (We would like to be able to say that $H_{2i+1}^{S^1}\cong \ZZ_k(e(1\otimes u^{k-1}))\gamma_{i-\frac{k}{2}(n-1)}\oplus \cdots \oplus \ZZ_3 e(1\otimes u^2)\gamma_{i-\frac{3}{2}(n-1)} \oplus \ZZ_2 e(1\otimes u)\gamma_{i-n+1}$, but there are extension issues that are difficult to resolve, so we cannot say which torsion group $H_{2i+1}^{S^1}$ should be.)\\
\indent Consider the $k+1$-th non-zero piece of the Gysin sequence:\\
\vspace{.2cm}
\begin{tikzpicture}
 [descr/.style={fill=white,inner sep=1.5pt}]
        \matrix (m) [
            matrix of math nodes,
            row sep=2em,
            column sep=.5em,
            text height=1.5ex, text depth=0.25ex
        ]
        {         
        & H^{S^1}_{(k+1)(n-1)+4} & H^{S^1}_{(k+1)(n-1)+2} & 0\\                  
             & H^{S^1}_{(k+1)(n-1)+3} & H^{S^1}_{(k+1)(n-1)+1}& 0\\     
        & H^{S^1}_{(k+1)(n-1)+2} & H^{S^1}_{(k+1)(n-1)} &  \ZZ (1\otimes u^k) \\     
        & H^{S^1}_{(k+1)(n-1)+1} & H^{S^1}_{(k+1)(n-1)-1}\cong \bigoplus_{j=2}^{k}\ZZ_j &\ZZ (a \otimes u^{k+1}) \\     
        & H^{S^1}_{(k+1)(n-1)} & H^{S^1}_{(k+1)(n-1)-2}\cong \ZZ(\gamma_{\frac{(k+1)(n-1)-2}{2}}) & 0\\    
        };

        \path[overlay,->, font=\scriptsize,>=latex]

        (m-1-2) edge  node[above] {$\cong$}(m-1-3)
        (m-1-3) edge[orange] node[above] {$M$}(m-1-4) 
        (m-1-4) edge[out=355,in=175,cyan] node[descr,yshift=0.3ex] {$e$} (m-2-2)
        (m-2-2) edge  (m-2-3)
        (m-2-3) edge[orange] node[above] {$M$}(m-2-4)
        (m-2-4) edge[out=355,in=175,cyan] node[descr,yshift=0.3ex] {$e$} (m-3-2)
        (m-3-2) edge (m-3-3)
        (m-3-3) edge[orange] node[above] {$M$} (m-3-4)
        (m-3-4) edge[out=355,in=175,cyan] node[descr,yshift=0.3ex] {$e$} (m-4-2)
        (m-4-2) edge (m-4-3)
        (m-4-3) edge[orange] node[above] {$M$} (m-4-4)
         (m-4-4) edge[out=355,in=175,cyan] node[descr,yshift=0.3ex] {$e$}(m-5-2)   
        (m-5-2) edge (m-5-3)
        (m-5-3) edge[orange] node[above] {$M$}(m-5-4)

    ;
\end{tikzpicture}\\
\normalfont
Thus, we can extract a short exact sequence from the last two lines of the Gysin sequence above, giving us $H^{S^1}_{(k+1)(n-1)}\cong \ZZ(e(a\otimes u^{k+1}))\oplus \ZZ (\gamma_{\frac{(k+1)(n-1)}{2}})$.
\normalfont
\indent It can be seen that $H^{S^1}_{(k+1)(n-1)+2}\cong \ZZ (\gamma_{\frac{(k+1)(n-1)+2}{2}})$ and that\\ $torsion(H^{S^1}_{(k+1)(n-1)+1})\cong t_{k+1}$. Since $H^{S^1}_{(k+1)(n-1)+1}$ is all torsion since it is sandwiched between a short exact sequence of torsion groups, we have $H^{S^1}_{(k+1)(n-1)+1}\cong  t_{k+1}$. Since loop homology $\mathbb{H}_{i}(LS^n)$ is zero for $(k+1)(n-1)+2-n \leq i \leq (k+2)(n-1)-1-n$), we obtain the analogous statements of~(\ref{eqn:oddinduction}) for $k+1$. In summary, we get the following theorem.
\thm For $n$ odd,
\begin{equation} 
H^{S^1}_{2i}(LS^n) \cong \left\{
\begin{array}{rl}
\ZZ \gamma_i & \text{if } (n-1)\nmid 2i,\\
\ZZ \gamma_i \oplus \ZZ(e(a\otimes u^i)) & \text{if } (n-1)\vert 2i. \nonumber
\end{array} \right.
\end{equation}
for $1\leq i \leq \frac{(k+1)(n-1)-2}{2}$ and 
\begin{equation}
H_{2i+1}^{S^1}(LS^n)\cong t_k  \nonumber
\end{equation}
for $\frac{k(n-1)}{2}\leq i \leq \frac{(k+1)(n-1)-2}{2}$, where $t_k$ is a torsion group of order $k!$. All other $j$ that does not fall into the above categories, we have that $H_i^{S^1}(LS^n) \cong H_{i-2}^{S^1}(LS^n)$
\normalfont
\subsection{The String Bracket for Odd Spheres}
~
\indent The string bracket is a degree $2-n$ map. The only possible non-zero bracket is of the generators $e(a\otimes u^i)$, since the marking map $M$ sends all other generators to zero.
\thm
\begin{align*}
[e(a\otimes u^i), e(a\otimes u^j)]=ije(1\otimes u^{i+j-2})
\end{align*}
where $e(1\otimes u^{i+j-2})$ is a generator of $\ZZ_{i+j-1}$, so the bracket is only zero when $i+j-1$ divides $ij$
\proof
\begin{align*}
[e(a\otimes u^i), e(a\otimes u^j)]=&e(M(e(a\otimes u^i))\bullet M(e(a\otimes u^j)) \nonumber \\
=&e(i(1\otimes u^{i-1})\bullet j(1\otimes u^{j-1})) \nonumber \\
=&e(ij(1\otimes u^{i+j-2})). \nonumber 
\end{align*}
\qed
\section{String Homology and the String Bracket of Even Spheres}
\subsection{Computations for $S^4$}
~
\indent Here we only use information from the Gysin sequence. As before, we know that
\begin{align*}
\mathbb{H}_*(LS^4,\ZZ)=\frac{\Lambda (b) \otimes \ZZ[a,v]}{(a^2,ab,2av)}
\end{align*}
where $|a|=-4$,$|b|=-1$, and $|v|=6$, so all of the generators look like $av^k$, $bv^k$, $v^k$, where $|av^k|=-4+6k$, $|bv^k|=-1+6k$, $|v^k|=6^k$ \cite{CJY}.\\
\indent We also know how the BV-operator acts,
\begin{align*}
\Delta (v^k)&=0\\
\Delta (av^k) &=0\\
\Delta (bv^k) &=(2k+1)v^k
\end{align*}
from \cite{Me}.
Let's consider the bottom of the Gysin sequence:\\
\begin{tikzpicture}[descr/.style={fill=white,inner sep=1.5pt}]
        \matrix (m) [
            matrix of math nodes,
            row sep=2em,
            column sep=2.5em,
            text height=1.5ex, text depth=0.25ex
        ]
        {      & H^{S^1}_7(LS^4) & H^{S^1}_5(LS^4) & \mathbb{H}_{2}(LS^4)\cong \ZZ_2 (av) \\
        & H^{S^1}_6(LS^4) & H^{S^1}_4(LS^4) & \mathbb{H}_{1}(LS^4)\cong  0 \\                
        & H^{S^1}_5(LS^4) & H^{S^1}_3(LS^4) & \mathbb{H}_{0}(LS^4)\cong \ZZ (1) \\
        & H^{S^1}_4(LS^4) & H^{S^1}_2(LS^4) & \mathbb{H}_{-1}(LS^4)\cong \ZZ (b) \\      
          & H^{S^1}_3(LS^4) & H^{S^1}_1(LS^4) & \mathbb{H}_{-2}(LS^4)\cong 0 \\          
            & H^{S^1}_2(LS^4) & H^{S^1}_0(LS^4) & \mathbb{H}_{-3}(LS^4)\cong 0 \\
              & H^{S^1}_1(LS^4) & 0 & \mathbb{H}_{-4}(LS^4)\cong \ZZ (a) \\
            & H^{S^1}_0(LS^4) &0  & 0. \\
        };

        \path[overlay,->, font=\scriptsize,>=latex]

        (m-1-2) edge (m-1-3)
        (m-1-3) edge[orange] node[above] {$M$}(m-1-4) 
        (m-1-4) edge[out=355,in=175,cyan] node[descr,yshift=0.3ex] {$e$} (m-2-2)
        (m-2-2) edge (m-2-3)
        (m-2-3) edge[orange] node[above] {$M$}(m-2-4)
        (m-2-4) edge[out=355,in=175,cyan] node[descr,yshift=0.3ex] {$e$} (m-3-2)
        (m-3-2) edge (m-3-3)
        (m-3-3) edge[orange] node[above] {$M$} (m-3-4)
        (m-3-4) edge[out=355,in=175,cyan] node[descr,yshift=0.3ex] {$e$} (m-4-2)
        (m-4-2) edge (m-4-3)
        (m-4-3) edge[orange] node[above] {$M$} (m-4-4)
         (m-4-4) edge[out=355,in=175,cyan] node[descr,yshift=0.3ex] {$e$}
(m-5-2)
        (m-5-2) edge (m-5-3)
        (m-5-3) edge[orange] node[above] {$M$} (m-5-4)
         (m-5-4) edge[out=355,in=175,cyan] node[descr,yshift=0.3ex] {$e$}                 
   (m-6-2)
        (m-6-2) edge (m-6-3)
        (m-6-3) edge[orange] node[above] {$M$} (m-6-4)
         (m-6-4) edge[out=355,in=175,cyan] node[descr,yshift=0.3ex] {$e$}        
 (m-7-2)
        (m-7-2) edge (m-7-3)
        (m-7-3) edge[orange] node[above] {$M$} (m-7-4)
         (m-7-4) edge[out=355,in=175,cyan] node[descr,yshift=0.3ex] {$e$}       
(m-8-2)   
        (m-8-2) edge (m-8-3)
        (m-8-3) edge[orange] node[above] {$M$}(m-8-4);
\end{tikzpicture}\\
The Gysin sequence, along with the BV-operator, allow us to determine that
\begin{align*}
H_0^{S^1}(LS^4)& \cong \ZZ (e(a))\\
H_1^{S^1}(LS^4)& \cong 0\\
H_2^{S^1}(LS^4)& \cong \ZZ (\gamma)\\
H_3^{S^1}(LS^4)& \cong \ZZ (e(b))\\
H_4^{S^1}(LS^4)& \cong \ZZ (\gamma_2)\\
H_5^{S^1}(LS^4)& \cong 0\\
H_6^{S^1}(LS^4)& \cong \ZZ (\gamma_3)\oplus \ZZ_2(e(av))\\
H_7^{S^1}(LS^4)& \cong 0.\\
\end{align*}
Here, $\gamma$ is the generator from $H_2(\CC P^{\infty})$.
Continuing up the Gysin sequence inductively, the $k$th piece of the sequence is as follows.\\
\begin{tikzpicture}[descr/.style={fill=white,inner sep=1.5pt}]
        \matrix (m) [
            matrix of math nodes,
            row sep=2em,
            column sep=1em,
            text height=1.5ex, text depth=0.25ex
        ]
        {      & H^{S^1}_{6k+7}(LS^4) & H^{S^1}_{6k+5}(LS^4) & \mathbb{H}_{6k+2}(LS^4)\cong \ZZ_2 (av^{k+1}) \\
        & H^{S^1}_{6k+6}(LS^4) & H^{S^1}_{6k+4}(LS^4) & \mathbb{H}_{6k+1}(LS^4)\cong  0 \\                
        & H^{S^1}_{6k+5}(LS^4) & H^{S^1}_{6k+3}(LS^4) & \mathbb{H}_{6k}(LS^4)\cong \ZZ (v^k) \\
        & H^{S^1}_{6k+4}(LS^4) & H^{S^1}_{6k+2}(LS^4) & \mathbb{H}_{6k-1}(LS^4)\cong \ZZ (bv^k) \\      
          & H^{S^1}_{6k+3}(LS^4) & H^{S^1}_{6k+1}(LS^4) & \mathbb{H}_{6k-2}(LS^4)\cong 0 \\          
            & H^{S^1}_{6k+2}(LS^4) & H^{S^1}_{6k}(LS^4) & \mathbb{H}_{6k-3}(LS^4)\cong 0 \\            
        };

        \path[overlay,->, font=\scriptsize,>=latex]

        (m-1-2) edge (m-1-3)
        (m-1-3) edge[orange] node[above] {$M$}(m-1-4) 
        (m-1-4) edge[out=355,in=175,cyan] node[descr,yshift=0.3ex] {$e$} (m-2-2)
        (m-2-2) edge (m-2-3)
        (m-2-3) edge[orange] node[above] {$M$}(m-2-4)
        (m-2-4) edge[out=355,in=175,cyan] node[descr,yshift=0.3ex] {$e$} (m-3-2)
        (m-3-2) edge (m-3-3)
        (m-3-3) edge[orange] node[above] {$M$} (m-3-4)
        (m-3-4) edge[out=355,in=175,cyan] node[descr,yshift=0.3ex] {$e$} (m-4-2)
        (m-4-2) edge (m-4-3)
        (m-4-3) edge[orange] node[above] {$M$} (m-4-4)
         (m-4-4) edge[out=355,in=175,cyan] node[descr,yshift=0.3ex] {$e$}
(m-5-2)
        (m-5-2) edge  (m-5-3)
        (m-5-3) edge[orange] node[above] {$M$} (m-5-4)
         (m-5-4) edge[out=355,in=175,cyan] node[descr,yshift=0.3ex] {$e$}                 
(m-6-2)
        (m-6-2) edge node[above] {$\cong$}(m-6-3)
        (m-6-3) edge[orange] node[above] {$M$} (m-6-4);
\end{tikzpicture}\\
Using the Poincar\'e polynomial for $H_*^{S^1}(LS^4,\ZZ_2)$ from \cite{We}, 
\begin{align*}
\left(\frac{1}{1-t^6}\right)\left( t^3+\frac{1+t^7}{1-t^2}\right)
\end{align*}
which we can rewrite in a more useful way as follows
\begin{align*}
\sum_{k=0}^{\infty}(k+1)(t^{6k}+t^{6k+2}+t^{6k+3}+t^{6k+4})+\sum_{k=1}^{\infty}k(t^{6k+1}+t^{6k+5}).
\end{align*}
From this, we see that $H_{6k+2}^{S^1}(LS^4,\ZZ_2)\cong \bigoplus\limits_{i=1}^{k} \ZZ_2$ Using the Universal Coefficient Theorem, 
\begin{align*}
0\rightarrow H_{6k+2}^{S^1}(LS^4,\ZZ) \otimes \ZZ_2 \rightarrow  \bigoplus_{i=1}^{k} \ZZ_2 \rightarrow Tor(H_{6k+1}^{S^1}(LS^4,\ZZ),\ZZ_2)\rightarrow 0
\end{align*}
Since $H_{6k+1}^{S^1}(LS^4),\ZZ)=0$, We have that $H_{6k+2}^{S^1}(LS^4,\ZZ)\otimes \ZZ_2 \cong \bigoplus\limits_{i=1}^{k} \ZZ_2$. We knot that $H_{6k+2}^{S^1}(LS^4,\ZZ)$ has a summand $\bigoplus_i \ZZ_{2^{l}}$ where $\sum l=k$, then we have
\begin{align*}
H_{6k+2}^{S^1}(LS^4,\ZZ)\otimes \ZZ_2 \cong\bigoplus_i \ZZ_{2^{l}} \cong \bigoplus_{i=1}^{k} \ZZ_2.
\end{align*}
So we must have that each $l=1$. Therefore, we have the following proposition.
\prop
\begin{align*}
H^{S^1}_{6k+7} & \cong 0\\
H^{S^1}_{6k+6} & \cong \ZZ (\gamma_{3k+3}) \oplus \ZZ_2({e(av^{k+1})})\bigoplus_{i=1}^{k+1} [\ZZ_2(e(av^{i})\gamma)] \oplus C_k\\
H^{S^1}_{6k+5} & \cong 0\\
H^{S^1}_{6k+4} & \cong \ZZ(\gamma_{3k+2}) \bigoplus_{i=1}^{k+1} [\ZZ_2(e(av^{i})\gamma)] \oplus C_k \\
H^{S^1}_{6k+3} & \cong \ZZ e(bv^k)\\
H^{S^1}_{6k+2} & \cong H^{S^{1}}_{6k} \cong \ZZ (\gamma_{3k})  \oplus C_{k-1}. \\
\end{align*}
$C_k$ is some torsion group of order $\prod_{i=1}^{k}(2i+1)$. The torsion comes from the fact that $M(e(bv^k))=\Delta(bv^k)=(2k+1)v^k$, which is where all of the odd torsion groups come from. 
\normalfont
\subsection{The String Bracket for $S^4$}
\indent The string bracket for $S^4$ is of degree $-2$.  We know that the marking map $M$ maps generators from $H_*(\CC P^{\infty})$ to zero (Lemma~\ref{lem:marking}), and it takes all the torsion to zero since for those cases, $M$ maps into zero or into a free group. Thus, the only possible case for the bracket to be nonzero is for the generators $e(bv^k)$. We have that
\begin{align*}
[e(bv^k),e(bv^l)]&=(-1)e(\Delta(bv^k)\bullet \Delta (bv^l))\\
&=-e((2k+1)v^k \bullet (2l+1)v^l)\\
&=-(2k+1)(2l+1)e(v^{k+l})\\
&=-4kle(v^{k+l}).
\end{align*}
\indent We know that $e(v^{k+l})$ has order $2(k+l)+1$,so the bracket is zero when $k \neq 0$, $l\neq 0$, and $(2k+2l+1)|(4kl+2k+2l+1)$, or when $4kl|(2k+2l+1)$, but the latter number is odd, so this can never happen. Thus, the bracket is always nontrivial in this case.\\
\indent When $k=0$ or $l=0$, the bracket is zero.
\subsection{Computations for $S^2$}
We just state results for the string homology computations for $S^2$, as the computations for $S^4$ are more illustrative. From \cite{CJY} , we have
\begin{align*}
\mathbb{H}_*(LS^2)=\frac{\Lambda(b)\otimes \ZZ [a,v]}{a^2, ab, 2av}
\end{align*}
where $|a|=-2$, $|b|=-1$, and $|v|=2$. The BV-operator acts as follows,
\begin{align*}
\Delta(v^k)&=0\\
\Delta(av^k)&=0\\
\Delta(bv^k)&=(2k+1)v^k+av^{k+1}.
\end{align*}
\prop
\begin{align*}
H_0^{S^1}(LS^2) & \cong \ZZ e(a)\\
H_2^{S^1}(LS^2) & \cong \ZZ e(a)\gamma \oplus \ZZ_2 e(av)\\
H_4^{S^1}(LS^2) & \cong \ZZ \oplus \ZZ_2 \oplus \ZZ_6 \\
H_{2i+1}^{S^1}(LS^2) & \cong \ZZ e(bv^i)\\
H_{2j}^{S^1}(LS^2) & \cong \ZZ \oplus C_k
\end{align*}
where $i \geq 0$, $j \geq 3$, and $C_k$ is a torsion group of order $\prod_{k=1}^{j-1}(4j-2-4k)$.
\normalfont
The string bracket is only non-zero in odd degrees.
\prop
\begin{align*}
[e(bv^i),e(bv^j)]=-4ije(v^{i+j})
\end{align*}
which is not always zero since $e(v^{i+j})$ is torsion, and all other brackets are zero.
\proof
\begin{align*}
[e(bv^i),e(bv^j)]& = (-1)^{2i+1-2}e(\Delta(bv^i)\bullet \Delta(bv^j))\\
&=-e(((2i+1)v^i+av^{i+1})\bullet((2j+1)v^j+av^{j+1}))\\
&=-e((2i+1)(2j+1)v^{i+j}+(2i+2j+2)av^{i+j+1})\\
&=-(2i+1)(2j+1)e(v^{i+j})\\
&=-4ije(v^{i+j}).
\end{align*}
\subsection{String Homology for Even Spheres}
We have from \cite{CJY} that, for $n$ even.
\begin{align*}
\mathbb{H}_*(LS^{n},\ZZ) \cong \frac{\Lambda(b) \otimes \ZZ[a,v]}{(a^2,ab,2av)}
\end{align*}
where $|a|=-n$, $|b|=-1$ and $|v|=2n-2$. By \cite{Me} we have
\begin{align*}
\Delta(v^k) &=0\\
\Delta(av^k) &=0\\
\Delta(bv^k) &= (2k+1)v^k.
\end{align*}
To keep track of things, $|av^k|=k(2n-2)-n$, $|bv^k|=k(2n-2)-1$, $|v^k|=k(2n-2)$. Let us consider the bottom of the Gysin sequence:\\

\begin{tikzpicture}[descr/.style={fill=white,inner sep=1.5pt}]
        \matrix (m) [
            matrix of math nodes,
            row sep=2em,
            column sep=2.5em,
            text height=1.5ex, text depth=0.25ex
        ]
        {       & H^{S^1}_{3n-1} & H^{S^1}_{3n-2} & \mathbb{H}_{2n-2}\cong \ZZ(v^2) \\                  
             & H^{S^1}_{3n-2} & H^{S^1}_{3n-4}& \mathbb{H}_{2n-3}\cong \ZZ (bv) \\   
             &\vdots &\vdots &\vdots\\  
        & H^{S^1}_{2n-1} & H^{S^1}_{2n-3}\cong 0 & \mathbb{H}_{n-2}\cong \ZZ_2 (av) \\     
        &\vdots &\vdots &\vdots\\
        & H^{S^1}_{n+1} & H^{S^1}_{n-1}\cong \ZZ & \mathbb{H}_{0}\cong \ZZ (v) \\     
        & H^{S^1}_{n} & H^{S^1}_{n-2} & \mathbb{H}_{-1}\cong \ZZ (b)\\        
        &\vdots &\vdots &\vdots\\
          & H^{S^1}_2 & H^{S^1}_0 & \mathbb{H}_{-n+1}\cong 0 \\
          & H^{S^1}_1 & 0 & \mathbb{H}_{-n}\cong \ZZ (a) \\
            & H^{S^1}_0 &0  &  \\
        };

        \path[overlay,->, font=\scriptsize,>=latex]

        (m-1-2) edge  node[above] {$\cong$}(m-1-3)
        (m-1-3) edge[orange] node[above] {$M$}(m-1-4) 
        (m-1-4) edge[out=355,in=175,cyan] node[descr,yshift=0.3ex] {$e$} (m-2-2)
        (m-2-2) edge node[above] {$\cong$} (m-2-3)
        (m-2-3) edge[orange] node[above] {$M$}(m-2-4)
        
        (m-4-2) edge (m-4-3)
        (m-4-3) edge[orange] node[above] {$M$} (m-4-4)
         
        (m-6-2) edge (m-6-3)
        (m-6-3) edge[orange] node[above] {$M$}(m-6-4)
         (m-6-4) edge[out=355,in=175,cyan] node[descr,yshift=0.3ex] {$e$}(m-7-2)   
        (m-7-2) edge  node[above] {$\cong$} (m-7-3)
        (m-7-3) edge[orange] node[above] {$M$}(m-7-4)
 (m-9-2)   
        (m-9-2) edge node[above] {$\cong$} (m-9-3)
        (m-9-3) edge[orange] node[above] {$M$}(m-9-4)
        (m-9-4) edge[out=355,in=175,cyan] node[descr,yshift=0.3ex] {$e$}(m-10-2)   
        (m-10-2) edge node[above] {$\cong$}(m-10-3)
        (m-10-3) edge[orange] node[above] {$M$}(m-10-4)
        (m-10-4) edge[out=355,in=175,cyan] node[descr,yshift=0.3ex] {$e,\cong$}(m-11-2)   
        (m-11-2) edge (m-11-3);
\end{tikzpicture}\\
from this sequence and knowledge of the BV-operator, we get
\pagebreak
\begin{align*}
H_0^{S^1}(LS^n)&\cong \ZZ a\\
H_1^{S^1}(LS^n)&\cong 0\\
H_2^{S^1}(LS^n)&\cong \ZZ \gamma \\
H_3^{S^1}(LS^n)&\cong 0 \\
H_4^{S^1}(LS^n)&\cong \ZZ \gamma_2\\
\vdots &\\
H_{n-1}^{S^1}(LS^n)&\cong \ZZ e(b)\\
H_{n}^{S^1}(LS^n)&\cong \ZZ e(v) \oplus \ZZ \gamma_{\frac{n}{2}} \\
\vdots &\\
H_{2n-2}^{S^1}(LS^n)&\cong \ZZ_2 e(av) \oplus \ZZ \oplus \ZZ \gamma_{\frac{2n-2}{2}} \\
H_{2n-1}^{S^1}(LS^n)&\cong 0 \\
H_{2n}^{S^1}(LS^n)&\cong \ZZ_2  \oplus \ZZ \oplus \ZZ \\
\vdots & \\
H_{3n-3}^{S^1}(LS^n)&\cong \ZZ e(bv)\\
H_{3n-2}^{S^1}(LS^n)&\cong \ZZ_2  \oplus \ZZ_3 e(v^3) \oplus \ZZ \oplus \ZZ \\
\vdots & 
\end{align*}
\begin{align*}
H_{4n-4}^{S^1}(LS^n)&\cong \ZZ_2 \oplus \ZZ_2 \oplus \ZZ_3 \oplus \ZZ \oplus \ZZ \\
H_{4n-3}^{S^1}(LS^n)&\cong 0\\
\vdots \\
H_{5n-5}^{S^1}(LS^n)&\cong \ZZ e(bv^2)\\
H_{5n-4}^{S^1}(LS^n)&\cong \ZZ_2 \oplus \ZZ_2 \oplus \ZZ_2 \oplus \ZZ_3 \oplus \ZZ \oplus \ZZ \\
\vdots &
\end{align*}
where all of the odd degree homology are isomorphic, and all even degree homology are isomorphic, or $H_i^{S^1}\cong H_{i-2}^{S^1}$m in the gaps denoted by the vertical dots. The $k$-th piece of the sequence is as follows:

\begin{tikzpicture}[descr/.style={fill=white,inner sep=1.5pt}]
        \matrix (m) [
            matrix of math nodes,
            row sep=2em,
            column sep=2.5em,
            text height=1.5ex, text depth=0.25ex
        ]
        {       & H^{S^1}_{k(2n-2)+n+1} & H^{S^1}_{k(2n-2)+n-1} & \mathbb{H}_{k(2n-2)}\cong \ZZ(v^k) \\                  
             & H^{S^1}_{k(2n-2)+n} & H^{S^1}_{k(2n-2)-2+n}& \mathbb{H}_{k(2n-2)-1}\cong \ZZ (bv^k) \\   
             &\vdots &\vdots &\vdots\\  
        & H^{S^1}_{k(2n-2)+1} & H^{S^1}_{k(2n-2)-1}\cong 0 & \mathbb{H}_{k(2n-2)-n}\cong \ZZ_2 (av^k) \\     
        };

        \path[overlay,->, font=\scriptsize,>=latex]

        (m-1-2) edge  node[above] {$\cong$}(m-1-3)
        (m-1-3) edge[orange] node[above] {$M$}(m-1-4) 
        (m-1-4) edge[out=355,in=175,cyan] node[descr,yshift=0.3ex] {$e$} (m-2-2)
        (m-2-2) edge node[above] {$\cong$} (m-2-3)
        (m-2-3) edge[orange] node[above] {$M$}(m-2-4)
        
        (m-4-2) edge (m-4-3)
        (m-4-3) edge[orange] node[above] {$M$} (m-4-4)
         
       ;
\end{tikzpicture}\\
inductively, we have that 
\begin{align*}
H_{k(2n-2)-2}^{S^1} \cong \ZZ_2^{k-1}\oplus C_k \oplus \ZZ \oplus \ZZ
\end{align*} 
where $C_k$ is a torsion group of order  $\prod_{i=1}^{k-1}(2i+1)$. The bottom of the above Gysin sequence gives the short exact sequence
\begin{align*}
0 \rightarrow \ZZ_2 av^k \rightarrow H^{S^1}_{k(2n-2)} \rightarrow H_{k(2n-2)-2}^{S^1} \rightarrow 0
\end{align*} 
which gives
\begin{align*}
H^{S^1}_{k(2n-2)}\cong \ZZ_2 e(av^k) \oplus \ZZ_2^{k-1}\oplus \left(\text{torsion group of order } \sum_{i=1}^{k-1}(2i+1)\right) \oplus \ZZ \oplus \ZZ.
\end{align*}
Note that the even torsion can be resolved using the results by Westerland in \cite{We} (as in the above example for $S^4$). From the top of the above Gysin sequence, we get the following.
\begin{align*}
H_{k(2n-2)+1}^{S^1}(LS^n)& \cong 0\\
\vdots & \\
H_{k(2n-2)+n-1}^{S^1}(LS^n) & \cong \ZZ e(bv^k)\\
H_{k(2n-2)+n}^{S^1}(LS^n) & \cong \ZZ_2^{k} \oplus C_k \oplus \ZZ \oplus \ZZ .
\end{align*}
Summarizing, we get the following theorem.
\thm Suppose $n$ is even. 
\begin{align*}
H_{k(2n-2)-2}^{S^1}(LS^n)&\cong \ZZ_2^{k-1}\oplus C_k \oplus \ZZ \oplus \ZZ\\
H_{k(2n-2)+1}^{S^1}(LS^n)& \cong 0\\
\vdots & \\
H_{k(2n-2)+n-1}^{S^1}(LS^n) & \cong \ZZ e(bv^k)\\
H_{k(2n-2)+n}^{S^1}(LS^n) & \cong \ZZ_2^{k} \oplus C_k \oplus \ZZ \oplus \ZZ .
\end{align*}
where all other unstated homology, we have $H_{i}^{S^1}(LS^n)\cong H_{i-1}^{S^1}(LS^n)$.
\normalfont.
\subsection{The String Bracket for Even Spheres}
\thm The string bracket is always zero except on the generators $e(bv^j)$,
\begin{align*}
[e(bv^k),e(bv^l)]=-(4kl+2k+2l+1)e(v^{k+l})
\end{align*}
but $e(v^{k+l})$ has order $2(k+l)+1$ so it is not always zero.
\proof 
\begin{align*}
[e(bv^k),e(bv^l)]&=(-1)^{k(2n-2)-1-n}e(M(e(bv^k))\bullet M(e(bv^l))\\
&=-(4kl+2k+2l+1)e(v^{k+l})
\end{align*}
\qed

\chapter{String Homology and String Bracket Computations for Surfaces}
~
\indent In this chapter we compute the string homology and string bracket for surfaces. 

\section{The Torus}
~
\indent We can compute the loop homology of the torus $T$ quite easily since we already know the loop homology of $S^1$. We compute the loop homology and BV-operator using the following,

\begin{align*}
\mathbb{H}_*(LT)\cong \mathbb{H_*(LS^1)}\otimes \mathbb{H_*(LS^1)}
\end{align*}
and that we can compute the BV-operator as follows:
\begin{align*}
\Delta_T(a\otimes b)=\Delta(a)\otimes b + (-1)^{|a|+1}a \otimes \Delta(b).
\end{align*}
We obtain the following:
\begin{align*}
\mathbb{H}_{-2}(LT)& \cong \bigoplus_{(n,m)\in \ZZ^2}\ZZ(1_{nm})\\
\mathbb{H}_{-1}(LT)& \cong \bigoplus_{(n,m)\in \ZZ^2}\ZZ x_{nm} \bigoplus_{(n,m)\in \ZZ^2} \ZZ y_{nm}\\
\mathbb{H}_{0}(LT)& \cong \bigoplus_{(n,m)\in \ZZ^2} \ZZ z_{nm}
\end{align*}
the loop product:
\begin{align*}
x_{nm}\bullet y_{kl}& =1_{n+k,m+l}\\
x_{nm}\bullet z_{kl}&=x_{n+k,m+l}\\
y_{nm}\bullet x_{kl}& = y_{n+k,m+l}\\
z_{nm}\bullet z_{kl} & = z_{n+k,m+l}
\end{align*}

and the BV-operator:

\begin{align*}
\Delta(1_{nm})& =nx_{nm}+my_{nm}\\
\Delta(x_{nm})& = nz_{nm}\\
\Delta(y_{nm}) & = -mz_{nm}\\
\Delta(z_{nm}) & =0.
\end{align*}

To calculate string homology of the torus, we use the Gysin sequence:

\begin{tikzpicture}[descr/.style={fill=white,inner sep=1.5pt}]
        \matrix (m) [
            matrix of math nodes,
            row sep=2em,
            column sep=.5em,
            text height=1.5ex, text depth=0.25ex
        ]
        { &  & &0 \\
        & H^{S^1}_3(LT)&H^{S^1}_1(LT) & \mathbb{H}_{0}(LT) \cong \bigoplus_{(n,m)\in \ZZ^2} \ZZ z_{nm}\\
          & H^{S^1}_2(LT) & H^{S^1}_0(LT) & \mathbb{H}_{-1}(LT) \cong \bigoplus_{(n,m)\in \ZZ^2}\ZZ x_{nm} \bigoplus_{(n,m)\in \ZZ^2} \ZZ y_{nm}\\
            & H^{S^1}_1(LT) & 0 & \mathbb{H}_{-2}(LT)\cong \bigoplus_{(n,m)\in \ZZ^1}\ZZ(1_{nm} \\
            & H^{S^1}_0(LT) & 0 & \\
           \\
        };

        \path[overlay,->, font=\scriptsize,>=latex]
    
        (m-1-4) edge[out=355,in=175,cyan] node[descr,yshift=0.3ex] {$e$} (m-2-2)
        (m-2-2) edge  node[above] {$c$} (m-2-3)
        (m-2-3) edge[orange] node[above] {$M$}(m-2-4)
        (m-2-4) edge[out=355,in=175,cyan] node[descr,yshift=0.3ex] {$e$} (m-3-2)
        (m-3-2) edge (m-3-3)
        (m-3-3) edge[orange] node[above] {$M$} (m-3-4)
        (m-3-4) edge[out=355,in=175,cyan] node[descr,yshift=0.3ex]        
        {$e$} (m-4-2)
        (m-4-2) edge (m-4-3)
        (m-4-3) edge[orange] node[above] {$M$} (m-4-4)
        (m-4-4) edge[out=355,in=175,cyan] node[descr,yshift=0.3ex]
 {$e$} (m-5-2)
        (m-5-2) edge (m-5-3)
         ;
\end{tikzpicture}\\
we get the following:

\begin{align*}
H_0^{S^1}(LT) \cong & \bigoplus_{(n,m)\in \ZZ^2}\ZZ e(1_{nm})\\
H_1^{S^1}(LT) \cong & \bigoplus_{\substack{(n,m)\in \ZZ^2 \\ gcd(n,m)=d \\ d\neq n,m}} \ZZ_d (q_{n_1}e(x_{nm})+ q_{m_1}e(y_{nm})) \bigoplus _{\substack{n|m\\ nd=m}}\ZZ_n (e(x_nm)+de(y_{nm})\\
& \bigoplus_{\substack{m|n \\ md=n}} \ZZ_m (de(x_{nm})+e(y_{nm})) \bigoplus_{(n,m)\in \Delta} \ZZ_n (e(x_{nn})+e(y_nn))\\
& \bigoplus_{gcd(n,m)=1}\ZZ(q_ne(x_{nm})+q_me(y_{nm})) \bigoplus_{gcd(n,m)=d}\ZZ(q_ne(x_{nm})+q_me(y_{nm}))\\
& \bigoplus_{n|m}\ZZ_e(y_{nm}) \bigoplus_{m|n}e(x_{nm})\bigoplus_{(n,n)\in \Delta}\ZZ e(y_{nn})\oplus \ZZ e(x_{00})\oplus e(y_{00})
\end{align*}
where for $H_1^{S^1}(LT)$, $(n,m)\in \ZZ-{(0,0)}$ for all $n,m$, and $\Delta$ is the diagonal in $\ZZ \times \ZZ$, and each $q_i$ are polynomials in the quotients that show up in the division algorithm for $n,m$. Note that the generator in the third to the last summand could also have been chosen to be $\ZZ e(x_nn)$. To calculate $H^{S^1}_2(LT)$, we look at the second and third line in the Gysin sequence above. Using that $\mathbb{H}_0(LT)$ is free, and $\Delta$ maps torsion to zero, the marking map is nonzero only on the free elements of $H_1^{S^1}(LT)$, we can extract the following short exact sequence:\\
\begin{align*}
0\rightarrow \frac{\mathbb{H}_0(LT)}{im(M)=im(\Delta)}\rightarrow H^{S^1}_2(LT) \rightarrow \ZZ 1_{00} \rightarrow 0
\end{align*}
where the last part of the sequence comes from the kernel of $M$ being generated by $1_{00}$. We get that
\begin{align*}
 \frac{\mathbb{H}_0(LT)}{im(M)}=\bigoplus_{(n,n)\in \Delta -(0,0)} \ZZ_n z_{nn} \bigoplus_{m|n} \ZZ_n z_{nm} \bigoplus_{n|m} \ZZ_m z_{nm} \bigoplus_{\substack{gcd(m,n)=d \\ d\neq n,m}} \ZZ_{q_n n -mq_m} z_{nm}. 
  \end{align*}
 Thus, we have
 \begin{align*}
 H_2^{S^1}(LT) \cong  \frac{\mathbb{H}_0(LT)}{im(M)} \oplus \ZZ 1_{00}\gamma
\end{align*}
From the very top of the Gysin sequence pictured above, we get that $H_3^{S^1}(LT) \cong ker(M)$, thus we have
\begin{align*}
H_3^{S^1}(LT) \cong tor(H_1^{S^1}(LT)) \oplus \ZZ x_{00}\gamma \oplus \ZZ y_{00}\gamma .
\end{align*}
Since loop homology higher than two is zero, we obtain the following isomorphisms:
\begin{align*}
H_{2k}^{S^1}(LT) & \cong H_2^{S^1}(LT)\\
H_{2k+1}^{S^1}(LT)& \cong H_3^{S^1}(LT)
\end{align*}
for all $k\geq 2$.
\subsection*{Torus Bracket Computations}
For $[-,-]: H_0^{S^1}(LT) \otimes H_0^{S^1}(LT) \rightarrow H_0^{S^1}(LT)$,
\begin{align*}
[e(1_{nm}),e(1_{kl})]=(nl-mk)e(1_{(n+k),(m+l)}).
\end{align*}
which corresponds to the Goldman bracket for the torus in Proposition~\ref{prop:torusstructure}. \\
For $[-,-]: H_0^{S^1}(LT) \otimes H_1^{S^1}(LT) \rightarrow H_1^{S^1}(LT)$,
\begin{align*}
[e(1_nm),e(q_1x_{nm}+q_2y_{kl})]=(q_1k-q_2l)e(nx_{n+k,m+l}+my_{n+k,m+l}).
\end{align*}
For $[-,-]: H_1^{S^1}(LT) \otimes H_1^{S^1}(LT) \rightarrow H_2^{S^1}(LT)$,
\begin{align*}
[e(q_1x_{nm}+q_2y_{nm}),e(q_3x_{kl}+q_4y_{kl})]=(q_1n-q_2m)(q_3k-q_4l)e(z_{n+k,m+l})
\end{align*}
which is torsion, and not always zero. All other brackets are zero.
\subsection*{String Lie Algebra Structure on the Torus}
~
\indent The center of the String Lie algebra can be directly computed from the above bracket results. Let $\mathfrak{g}(T)$ denote the String Lie algebra for the closed torus, the center of the Lie algebra is as follows:
\begin{align*}
Z(\mathfrak{g}(T))& \cong H_*^(S^1)(L_0(T) \bigoplus_{k\in \ZZ - \{ 0 \} } \ZZ e(y_{k0})\bigoplus_{l\in \ZZ - \{0\}} \ZZ e(x_{0l}) \\ & \bigoplus_{(m,n)\in \ZZ ^2 -\{(0,0)\}} \ZZ e(z_{mn}) \oplus \text{ tor}(H_1^{S^1}(LT)) \oplus H_2^{S^1}(LT)
\end{align*}
where $L_0(T)$ denotes the connected component of the loop space $LT$ containing the contractible loops.
\section{Genus $g>1$}
~
\indent Using that $L(BG)=\amalg_{[\gamma]} BC(\gamma)$, for a closed, orientable surface $\Sigma_g$ of genus $g>1$, we have
\begin{align*}
L\Sigma_g=LB(\pi_1(\Sigma_g,*))=\coprod_{[\gamma]\in \hat{\pi}} BC(\gamma)
\end{align*}
where $[\gamma ]$ is a conjugacy class in $\pi_1(\Sigma_g)$, or an element of $\hat{\pi}$. By Kupers \cite{Ku}, $C(\gamma)\cong \ZZ$ for $\gamma \neq e$, where $e$ is the identity in $\pi_1(\Sigma_g)$. This gives us that $BC(\gamma)=S^1$. For the centralizer of $e$, $C(e)=\pi_1(\Sigma_g)$ since everything commutes with $e$, so $BC(e)=\Sigma_g$. Therefore, we have
\begin{align*}
L\Sigma_g=\coprod_{[\gamma]\neq e \in \hat{\pi}}S^1 \cup \Sigma_g.
\end{align*}
We can follow \cite{Ku}, and also knowing the loop homology of $S^1$ and $\Sigma_g$, we get the following 
\thm (Kupers, \cite{Ku}, Theorem 2.2)
\begin{align*}
\mathbb{H}_{-2}(L\Sigma _g) &= \bigoplus_{\gamma \in \hat{\pi}} \ZZ [\gamma ]\\
\mathbb{H}_{-1}(L\Sigma _g) &=H_1(\Sigma _g)\bigoplus_{\gamma \in \hat{\pi}-\{e\}} \ZZ \beta_{\gamma}\\
\mathbb{H}_0(L\Sigma _g) &=\ZZ 1\\
\mathbb{H}_i(L \Sigma_g) &=0 \hspace*{.5cm}\text{for}\hspace*{.2cm} i>0.
\end{align*}
\normalfont
Let $H_1(\Sigma_g)=\bigoplus\limits_{i=1}^g\ZZ a_i \bigoplus\limits_{j=1}^g\ZZ b_i$. Let $\kappa_{\gamma}$ be the generator of $C[\gamma ]$, then $\gamma = \kappa_{\gamma}^{l_{\gamma}}$, where $l_{\gamma}\in \ZZ$.  The loop product can be computed as follows \cite{Ku},
\begin{align}
1 \bullet 1 &=1\\
\beta_{\gamma} \bullet 1 &= \beta_{\gamma}\\
a_i \bullet 1 &= a_i\\
b_j \bullet 1 &= b_i\\
a_i \bullet a_j &=0\\
b_i \bullet b_j &=0\\
a_i \bullet b_j &= \delta_{ij}\\
[\gamma] \bullet 1 &= [\gamma]\\
\beta_{[\gamma_1]}\bullet \beta_{[\gamma_2]} &= \frac{[\beta_{[\gamma_1]}, \beta_{[\gamma_2]}]}{l_{\gamma_1}\cdot l_{\gamma_2}}
\end{align}
where (4.3)-(4.6) are just the homology intersection product on $H_1(\Sigma_g)$ and the bracket in (4.9) is the Goldman bracket. We also can compute the BV-operator as in \cite{Ku}, where the only non-trivial result is $\Delta ([\gamma])=l_{\gamma}\beta_{\gamma}$.\\
\indent As before, we can use the Gysin sequence to compute string homology. Consider the bottom of the Gysin sequence,\\
\begin{tikzpicture}[descr/.style={fill=white,inner sep=1.5pt}]
        \matrix (m) [
            matrix of math nodes,
            row sep=2em,
            column sep=1em,
            text height=1.5ex, text depth=0.25ex
        ]
        { &  & &0 \\
        & H^{S^1}_3(L\Sigma_g)&H^{S^1}_1(L\Sigma_g) & \mathbb{H}_{0}(L\Sigma_g) \cong \ZZ 1 \\
          & H^{S^1}_2(L\Sigma_g) & H^{S^1}_0(L\Sigma_g) & \mathbb{H}_{-1}(L\Sigma_g)\cong \bigoplus\limits_{\gamma \in \hat{\pi}-\{e\}} \ZZ \beta_{\gamma} \oplus H_1(L\Sigma_g)\\
            & H^{S^1}_1(L\Sigma_g) & 0 & \mathbb{H}_{-2}(L\Sigma_g)\cong \bigoplus\limits_{\gamma \in \hat{\pi}} \ZZ [\gamma]\\
            & H^{S^1}_0(L\Sigma_g) & 0 & \\
           \\
        };

        \path[overlay,->, font=\scriptsize,>=latex]
    
        (m-1-4) edge[out=345,in=160,cyan] node[descr,yshift=0.1ex] {$e$} (m-2-2)
        (m-2-2) edge  node[above] {$c$} (m-2-3)
        (m-2-3) edge[orange] node[above] {$M$}(m-2-4)
        (m-2-4) edge[out=345,in=160,cyan] node[descr,yshift=0.3ex] {$e$} (m-3-2)
        (m-3-2) edge node[above] {$c$}(m-3-3) 
        (m-3-3) edge[orange] node[above] {$M$} (m-3-4)
        (m-3-4) edge[out=345,in=165,cyan] node[descr,yshift=0.3ex]        
        {$e$} (m-4-2)
        (m-4-2) edge (m-4-3)
        (m-4-3) edge[orange] node[above] {$M$} (m-4-4)
        (m-4-4) edge[out=345,in=160,cyan] node[descr,yshift=0.3ex]
 {$e$} (m-5-2)
        (m-5-2) edge (m-5-3)
         ;
\end{tikzpicture}\\
we get that $H_0^{S^1}(L\Sigma_g)\cong \bigoplus\limits_{\gamma \in \hat{\pi}} \ZZ e([\gamma])$. Since $M \circ e ([\gamma])=\Delta ([\gamma])=l_{\gamma}\beta_{[\gamma]}$ for $[\gamma]\neq e$, $\Delta ([e])=0$ and $\Delta (a_i)=\Delta(b_j)=0$, we get that $H_1^{S^1}(L\Sigma_g) \cong \mathbb{H}_{-1}(L \Sigma_g)/ \bigoplus\limits_{\gamma \in \hat{\pi}-\{e\}}(\ZZ l_{\gamma}\beta_{\gamma})$. Since $M: H_1^{S^1}(L \Sigma_g) \rightarrow \ZZ 1$ is the zero map, $H_2^{S^1}(L \Sigma_g)$ sits in the following short exact sequence
\begin{align*}
0 \rightarrow \ZZ_1 \rightarrow H_2^{S^1}(L\Sigma_g) \rightarrow im(c) \rightarrow 0.
\end{align*}
Since $im(c)=ker(M)=\ZZ e([e])$, then the above exact sequence splits and $H_2^{S^1}(L\Sigma_g)\cong \ZZ e(1) \oplus \ZZ s$ (EXPLAIN $s$). Therefore, we have the following:
\prop \begin{align*}
H_0^{S^1}(L\Sigma_g)& \cong \bigoplus\limits_{\gamma \in \hat{\pi}} \ZZ e([\gamma])\\
H_1^{S^1}(L\Sigma_g) & \cong \bigoplus\limits_{i=1}^g\ZZ e(a_i) \bigoplus\limits_{j=1}^g\ZZ e(b_i) \bigoplus\limits_{\gamma \in \hat{\pi}-\{e\}}\ZZ_{l_{\gamma}}e(\beta_{\gamma})\\
H_2^{S^1}(L\Sigma_g) & \cong \ZZ e(1) \oplus \ZZ s\\
H_{2i+1}^{S^1}(L \Sigma_g) & \cong H_1^{S^1}(L \Sigma_g), \hspace*{.2cm} i\geq 1\\
H_{2i}^{S^1}(L \Sigma_g) & \cong H_2^{S^1}(L \Sigma_g), \hspace*{.2cm} i \geq 1.
\end{align*}
\normalfont
From this, we can compute the string bracket for $\Sigma_g$. The only non-trivial bracket is the Goldman bracket.
\prop
\begin{align}\label{algn:surfacegoldman}
[e([\gamma_1]),e([\gamma_2])]=[\gamma_1,\gamma_2] 
\end{align}
\normalfont
The second bracket in \ref{algn:surfacegoldman} is the Goldman bracket. So the only non-trivial string bracket of string homology of $\Sigma_g$ is the Goldman bracket as in Chapter 1.

\begin{singlespace}

\end{singlespace}
\end{document}